\documentclass[12pt,reqno,a4wide]{article}
\usepackage{a4wide}
\usepackage{eurosym}
\usepackage{comment}
\usepackage{mathptmx}
\usepackage{amssymb}
\usepackage{amsfonts}
\usepackage{amsthm}
\usepackage{amsmath}
\usepackage{dsfont}
\usepackage{graphicx}
\usepackage{caption}
\usepackage{shadow}
\usepackage[all]{xy}
\usepackage{enumerate}
\usepackage{makeidx}
\usepackage{mathrsfs}
\usepackage{upgreek}
\usepackage{stmaryrd}
\usepackage[utf8]{inputenc}
\usepackage{array}
\usepackage{tikz-cd}
\usepackage{thumbpdf}
\usepackage[pagebackref]{hyperref}
\usepackage{textcomp}
\usepackage{chapterbib}
\usepackage{graphics}
\usepackage{longtable}
\usepackage{epsf}
\usepackage{bm}
\usepackage{adjustbox}
\usepackage{centernot}
\usepackage{setspace}
\usepackage{imakeidx}
\usepackage{tikz}
\usepackage{tikz-cd}
\usepackage{float}

\makeindex
\newtheorem{theorem}{Theorem}[section]
\newtheorem{lem}[theorem]{Lemma}
\newtheorem{thm}[theorem]{Theorem}
\newtheorem{prop}[theorem]{Proposition}
\newtheorem{cor}[theorem]{Corollary}
\theoremstyle{definition}
\newtheorem*{Beweis}{Proof}
\newtheorem*{notatation}{Notation}
\newtheorem{defn}[theorem]{Definition}

\newtheorem{defns}[theorem]{Definitions}
\newtheorem{rem}[theorem]{Remark}

\newtheorem{punto}[theorem]{}
\newtheorem*{notation}{Notation}

\theoremstyle{remark}

\newtheorem{ex}[theorem]{Example}
\newtheorem{exs}[theorem]{Examples}
\newtheorem{remarks}[theorem]{Remarks}
\input{tcilatex}
\begin{document}

\title{Lower Separation Axioms for $X$-top Lattices\thanks{%
MSC2010: Primary 06A15, 6B30; Secondary: 13C13, 16D10, 16Y60 } \thanks{%
Keywords: Lattices, $X$-Top Lattice, Zariski Topology, Separation Axioms,
Quarter Separation Axioms, Zero Dimensional (Semi)Rings.} {\small \thanks{
This is an extended version of results extracted from the M.Sc. thesis of
the second author under the supervision of Prof. J. Abuhlail.}}}
\author{$%
\begin{array}{cc}
\text{Jawad Abuhlail}{\small \thanks{\text{Corresponding Author}}} & \text{%
Abdulmushin Alfaraj} \\
\text{{\small abuhlail@kfupm.edu.sa; abuhlail@yahoo.de}} & \text{%
abdulmuhsinalfaraj@hotmail.com} \\
\text{{\small Department of Mathematics}} & \text{Department of Mathematical
Sciences} \\
\text{{\small King Fahd University of Petroleum $\&$ Minerals}} & \text{%
University of Bath} \\
\text{{\small 31261 Dhahran, KSA}} & \text{Celeverton Down, Bath, BA27AY, UK}%
\end{array}%
$}
\date{\today }
\maketitle

\begin{abstract}
We study separation axioms for $X$-top-lattices (i.e. lattices $L$ for which
a given subset $X\subseteq L\backslash \{1\}$ admits a \emph{Zariski-like
topology}). Such spaces are $T_{0}$ and usually far away from being $T_{2}.$%
We give graphical characterizations for an $X$-top-lattice to be $T_{1},$ $%
T_{\frac{1}{4}},$ $T_{\frac{1}{2}},$ $T_{\frac{3}{4}}$ and provide several
families of examples/counterexamples that illustrate our results. We apply
our results mainly to the prime (resp. maximal, minimal) spectra of prime
(resp. maximal, minimal) ideals of commutative (semi)rings.
\end{abstract}

\section*{Introduction}

The prime spectrum $Spec(R)$ of a commutative ring $R$ attains the so called
\emph{Zariski topology}. This topology is very important in Commutative
Algebra (e.g. \cite{AM1969}) as it provides access to topological tools in
addition to the conventional algebraic ones. Moreover, this topology is
heavily used in Algebraic Geometry (e.g. \cite{Eis1995}).

Several authors attempted to generalize the standard Zariski topology to the
spectra of special classes of submodules of a module over a (commutative)\
ring. That enriched the Theory of Modules with many new classes of (\emph{co}%
)\emph{prime} submodules (e.g. \cite{Wij2006}, \cite{Hro2016}). The main
problem was that the closed varieties in such spectra are \emph{not
necessarily closed under finite unions}. Several authors tried to find
conditions under which such spectra have a Zariski-like topology. The first
attempts were by McCasland Smith (\cite{MMS1997}, \cite{MS2006}), who
introduced the notion of \emph{top modules }(see also Lu \cite{Lu1999}).
Moreover, several authors, including Abuhlail (e.g. \cite{Abu2011}, \cite%
{Abu2011-CA}, \cite{Abu2015}), studied \emph{Zariski-like topologies} for
different classes of \emph{(co)prime submodules} of a module over a ring
(e.g. \cite{BH2008}, \cite{Tek2009}). These investigations were extended by
several authors to special classes of \emph{(co)prime subsemimodules} of a
given semimodule over a semiring (e.g. \cite{AUT2013}, \cite{HPH2021}).

This led Abuhlail and Lomp \cite{AL2016} (see also \cite{AL2013}) to
consider the more general framework of what they called an $X$\emph{-top
lattice} (a \emph{complete lattice} $L$ for which an appropriate subset $%
X\subseteq L\backslash \{1\}$ attains a Zariski-like topology). Applications
are, in particular, to the lattice $Sub(_{R}M)$ of sub(semi)modules of a
given left (semi)module $M$ over a (semi)ring $R$. For the special case $M=R$%
, a (semi)ring, and $X=Spec(R),$ one recovers the Zariski topology on the
spectrum of prime ideals of $R.$

Abuhlail and Hroub \cite{AH2019} (see also \cite{Hro2016}) studied further
the properties of $X$-top lattices. They proved several results in this
general framework that recovered many results by different authors in
special contexts. Moreover, they studied the interplay between the algebraic
properties of the module/ring involved and the topological properties of
these Zariski-like topology as well as conditions under which such a
topology is \emph{spectral} (in the sense of Hochster \cite{Hoc1969}).

Our goal in this paper is to initiate a \emph{comprehensive study} of the
separation axioms of $X$-top lattices. Such an $X$ is $T_{1}$ if and only if
$X$ is \emph{Krull zero-dimensional} (i.e. $Max(X)=Min(X)$). Since such
topological spaces are $T_{0}$, we focus on separation axioms between $T_{0}$
and $T_{1},$ namely the \emph{quarter separation axioms}: $T_{\frac{1}{4}}$,
$T_{\frac{1}{2}}$ and $T_{\frac{3}{4}}$. Applications could be given to
several special classes of (semi)modules over (semi)rings. However, we
restrict ourselves to (semi)rings.

In Section 1, we recall general preliminaries from the theories of \emph{%
lattices}, \emph{semirings}, $X$\emph{-top lattices} and \emph{general
topology}. In Section 2, we consider the $X$-top lattices for which $X$ is $%
T_{1}.$ We demonstrate in Proposition \ref{dim-0} that $X$ is $T_{1}$ if and
only if $X$ is $K.\dim (X)=0$. A key result is Proposition \ref{iso}, which
characterizes the isolated points in $X,$ and used in Theorem \ref{discrete}
to characterize the discreteness of $X.$ We introduce several classes of
(semi)rings $R$ and characterize them through the discreteness of special
subspaces of $Spec(R)$: we call a minimal prime (a maximal ideal) $P$ \emph{%
absolutely minimal} (\emph{barely maximal}) iff $P$ does \emph{not} contain
the intersection of the \emph{other} minimal primes (maximal ideals). We
call a (semi)ring $R$ an \textbf{AMin (semi)ring} (a \textbf{BMax (semi)ring}%
) iff every minimal prime (maximal ideal) of $R$ is absolutely minimal
(barely maximal). In Corollary \ref{semilocal-discrete}, we deduce that a
commutative (semi)ring $R$ is a BMax-(semi)ring if and only if $Max(R)$ is
discrete and prove that this is equivalent to $R$ being an \textbf{%
FMax-(semi)ring}, i.e. a (semi)ring with finitely many maximal ideals. We
deduce in Corollary \ref{Min-discrete} that a commutative (semi)ring $R$ is
an AMin-(semi)ring if and only if $Min(R)$ is discrete. While an
FMin-(semi)ring) (\emph{i.e.} one with finitely many minimal prime ideals)
is clearly an AMin-(semi)ring, Example \ref{AMin-Not-FMin} shows that the
converse is \emph{not} true. In Theorem \ref{PAM} which characterizes a
commutative (semi)ring with discrete prime spectrum as a Krull $0$%
-dimensional AMin-(semi)ring (BMax-(semi)ring). In Theorem \ref{Artinian},
we obtain a new characterization of a commutative Artinian ring as a
Noetherian in which every prime ideal is absolutely minimal (barely
maximal). Moreover, we provide a series of examples showing that some
properties of special commutative semirings may fail for \emph{proper
semirings} of the same type. While a commutative ring $R$ is $\pi $-regular
if and only if $K.\dim (R)=0,$ we introduce in Example \ref{S3} a finite
(Artinian/Noetherian) subtractive local reduced idempotent (whence von
Neumann regular, $\pi $-regular)\ commutative semidomain with \emph{non-zero
Krull dimension}. For a commutative ring $R$ we have:
\begin{equation*}
R\text{ is Noetherian }\&\text{ }\pi \text{-regular}\Longrightarrow R\text{
is Artinian}\Longrightarrow R\text{ is Noetherian }\&\text{ }K.\dim (R)=0
\end{equation*}%
(and conversely, cf. Theorem \ref{Artinian}). Both implications fail even
for \emph{subtractive local} semidomains as Examples \ref{Noeth-not-Art-dim0}
and \ref{Art-not-Noeth} demonstrate.

In Section 3, we provide different characterizations of $X$-top lattices for
which $X$ satisfies one of the quarter separation axioms studied for
semirings in \cite{PRV2009}. A topological space $X$ is said to be $T_{\frac{%
1}{4}}$ (resp. $T_{\frac{1}{2}},$ $T_{\frac{3}{4}}$)\ iff every $x\in X$ is
either closed or \emph{kerneled} (resp. either closed or \emph{isolated},
either closed or \emph{regular open}). Let $L$ be an $X$-top lattice for
some $X\subseteq L\backslash \{1\}.$ Theorem \ref{t1/4t} characterizes the $%
X $-top lattices for which $X$ is $T_{\frac{1}{4}}$ as those with $K.\dim
(X)\leq 1.$ Theorem \ref{T1/2} characterizes the $X$-top lattices for which $%
X$ is $T_{\frac{1}{2}}$ as those for which every element of $X$ is either
maximal or \emph{simultaneously} minimal $\&$\emph{\ }completely strongly $X$%
-irreducible. Theorem \ref{3/4T} characterizes the $X$-top lattices for
which $X$ is $T_{\frac{3}{4}}$ as those for which each element of $X$ is
either maximal or \emph{simultaneously} minimal, completely strongly $X$%
-irreducible and excluded. Several sufficient/necessary conditions for such
an $X$ to satisfy the quarter separation axioms are given in case $X$ is a
finite forest of \emph{trees} $\mathcal{T}_{n}$ with finite base and/or
\emph{dual trees} $\mathcal{V}_{m}$ with finite covers (e.g. Theorem \ref%
{t3/4t}). This enables us to construct several families of
examples/counterexamples that illustrate our results. We end this section
with Corollary \ref{Bni} which provides applications to the prime spectra of
the semirings $B(n,i)$ introduced by F. Alarcon and D. D. Anderson \cite%
{AA1994} (as a generalization of the rings $\mathbb{Z}_{n}$).

\section{Preliminaries}

In this section, we introduce some relevant definitions and results from
theories of lattices, semirings, $X$-top lattices and general topology.

\subsection*{Lattices}

We recall the relevant definitions and results from Lattice Theory. We
follow \cite{Gra2010} (unless otherwise stated explicitly).

\begin{punto}
A \textbf{lattice} is a \emph{partially ordered set} $(L,\leq )$ such that
the \emph{infimum} and the \emph{supremum} of any two elements exist. For a
lattice $(L,\leq ),$ we set%
\begin{equation*}
a\vee b:=\sup \{a,b\}\text{ and }a\wedge b:=\inf \{a,b\}\text{ }\forall \
a,b\in L.
\end{equation*}%
We call $(L,\wedge )$ a \textbf{meet-semilattice} of $L$ and $(L,\vee )$ a
\textbf{join-semilattice} of $L.$ For $x,y\in X,$ we set ${y||x}$ iff $x$
and $y$ are \textbf{not comparable} (i.e. $x\nleqq y$ and $y\nleqq x$).
\end{punto}

\begin{defns}
A lattice $\mathcal{L}=(L,\wedge ,\vee )$ is said to be a

\begin{enumerate}
\item \textbf{bounded lattice }iff $L$ has a \emph{largest element} $1$ and
a \emph{smallest element} $0$.

\item \textbf{complete lattice} iff arbitrary meets and joins of elements of
$L$ exist in $L$ (in this case, $1=\bigvee\limits_{x\in L}x$ and $%
0=\bigwedge\limits_{x\in L}x,$ i.e. every complete lattice is bounded).

\item \textbf{distributive lattice} iff $\forall \ x,y,z\in L$:%
\begin{equation*}
x\wedge (y\vee z)=(x\wedge y)\vee (x\wedge z)\text{ (equivalently }x\vee
(y\wedge z)=(x\vee y)\wedge (x\vee z)\text{).}
\end{equation*}
\end{enumerate}
\end{defns}

\begin{punto}
Let $(L,\wedge )$ be a \emph{complete} meet-semilattice and $B,C\subseteq L.$
We call $q\in C:$

\textbf{(completely) }$B$-\textbf{irreducible in }$L$ iff for any $A\underset%
{\text{finite}}{\subseteq }B$ ($A\subseteq B$) with $q\leq a$ for all $a\in
A,$ we have:\ $\bigwedge\limits_{a\in A}a=q\Longrightarrow a=q$ for some $%
a\in A;$

\textbf{(completely) strongly }$B$-\textbf{irreducible in }$L$ iff for any $A%
\underset{\text{finite}}{\subseteq }B$ ($A\subseteq B$), we have: $%
\bigwedge\limits_{a\in A}a\leq q\Longrightarrow a\leq q$ for some $a\in A.$

With $I^{B}(C)$ (resp. $CI^{B}(C),$ $SI^{B}(C),$ $CSI^{B}(C)$), we denote
the set of $B$-irreducible (resp. completely $B$-irreducible, strongly $B$%
-irreducible, completely strongly $B$-irreducible) elements of $C.$ We drop
the superscript $B$ if it's clear from the context.
\end{punto}

\begin{punto}
Let $\mathcal{L}=(L,\wedge ,\vee ,1,0)$ be a \emph{complete} lattice and $%
X\subseteq L.$ We call $q\in X:$

\textbf{absolutely }$X$\textbf{-minimal} iff $q\in Min(X)$ and $%
\bigwedge\limits_{Min(X)\backslash \{q\}}a\nleqq q;$

\textbf{barely }$X$\textbf{-maximal} iff $q\in Max(X)$ and $%
\bigwedge\limits_{Max(X)\backslash \{q\}}a\nleqq q.$

With $AMin(X)=CSI^{Min(X)}(Min(X))$ (resp. $BMax(X)=CIS^{Max(X)}(Max(X))$),
we denote the set of absolutely $X$-minimal (barely $X$-maximal) elements of
$X.$

We say that $X$ is \textbf{AMin} (\textbf{BMax}) iff $Min(X)=AMin(X)\ $($%
Max(X)=BMax(X)$). We say that $X$ is \textbf{PAMin} (\textbf{PBMax}) iff $%
X=AMin(X$ ($X=BMax$).
\end{punto}

\begin{punto}
Let $\mathcal{L}=(L,\wedge ,\vee ,1,0)$ be a \emph{complete} lattice. We say
that

$L$ has the \textbf{max-property} iff $\underset{A\underset{\text{finite}}{%
\subseteq }Max(L)\backslash \{q\}}{\bigwedge a}\nleq q$ for every $q\in
Max(L),$ equivalently, $Max(L)=SI^{Max(L)}(Max(L));$

$L$ has the \textbf{complete max-property} iff $\underset{Max(L)\backslash
\{q\}}{\bigwedge a}\nleq q$ for every $q\in Max(L),$ equivalently, $%
Max(L)=CSI^{Max(L)}(Max(L)).$
\end{punto}

\begin{punto}
Let $\mathcal{L}=(L,\wedge ,\vee ,1,0)$ be a \emph{bounded} lattice and $%
X\subseteq A\subseteq L.$ We say that $A$ is

\begin{enumerate}
\item $X$-\textbf{atomic} iff for every $a\in A:$ there exists $m\in Min(X)$
such that $m\leq a;$

\item $X$-\textbf{coatomic} iff for every $a\in A:$ there exists $\mathfrak{m%
}\in Max(X)$ such that $a\leq \mathfrak{m}.$

We say that $X\subseteq L$ is \textbf{atomic} (\textbf{coatomic}) iff $X$ is
$X$-atomic ($X$-coatomic).
\end{enumerate}
\end{punto}

\subsection*{X-top Lattices}


We recall the some definitions and notation from the \emph{Theory }$X$\emph{%
-Top Lattices.} We follow \cite{AL2016} (unless otherwise stated explicitly).

\begin{notation}
Let $\mathcal{L}=(L,\wedge ,\vee ,1,0)$ be a complete lattice and $%
X\subseteq L\backslash \{1\}$. For any $a\in L$, we define%
\begin{equation*}
\begin{tabular}{lllllll}
$V_{X}(a)$ & $:=$ & $\{x\in X\ |\ a\leq x\};$ &  & $D_{X}(a)$ & $:=$ & $%
X\backslash V_{X}(a).$%
\end{tabular}%
\end{equation*}%
We call $V_{X}(a)$ the \emph{variety} of $a$ in $X.$ Moreover, we set
\begin{equation*}
V_{X}(L):=\{V_{X}(a)\ |\ a\in L\}\text{ and }\tau _{X}(L)=\{D_{X}(a)\mid
a\in L\}.
\end{equation*}%
We drop the subscript $X$ if it's clear from the context.
\end{notation}

\begin{punto}
\label{IV-def}Let $\mathcal{L}=(L,\wedge ,\vee ,1,0)$ be a \emph{complete
lattice} and $X\subseteq L\backslash \{1\}.$ Notice that $V_{X}(0)=X,$ $%
V_{X}(1)=\emptyset $ and $V_{X}(L)$ is closed under arbitrary intersections
as $\bigcap\limits_{a\in A}(V_{X}(a))=V_{X}(\bigvee\limits_{a\in A}a)$ for
any $A\subseteq L.$ We say that $\mathcal{L}$ is an $X$\textbf{-top lattice}
\cite{AL2016} iff $V_{X}(L)$ is closed under finite unions (equivalently, $%
\tau _{X}(L)$ is closed under finite intersections). In this case, $(X,\tau
_{X}(L))$ is a topological space, which we call a \textbf{Zariski-like
topology}). \newline
Consider $X\subseteq L\backslash \{1\}$. For any $Y\subseteq X$ and $a\in L,$
we set
\begin{equation*}
I_{X}(Y):=\bigwedge_{y\in Y}y\text{ and }\sqrt[X]{a}:=I_{X}(V_{X}(a)).
\end{equation*}%
It can be shown that $\overline{Y}=V_{X}(I_{X}(Y))$ (e.g. \cite[Lemma 1.8]%
{AH2019}. We say that $a\in L$ is $X$-\textbf{radical} iff $\sqrt[X]{a}=a$
and denote the \textbf{set of }$X$\textbf{-radical elements} of $L$ by $%
C^{X}(L)$, i.e.%
\begin{equation*}
C^{X}(L):=\{a\in L\text{ }|\text{ }a=\bigwedge_{y\in V_{X}(a)}y\}.
\end{equation*}%
Clearly, $X\subseteq C^{X}(L)$ and $(C^{X}(L),\wedge )$ is a
meet-semilattice. We drop the superscript $X$ if it's clear from the context.
\end{punto}

We recall a fundamental characterization of $X$-top lattices by Abuhlail and
Lomp \cite{AL2016}.

\begin{thm}
\label{xct}\emph{(\cite[\emph{Theorem 2.2}]{AL2016})} Let $\mathcal{L}%
=(L,\wedge ,\vee ,1,0)$ be a complete lattice and $X\subseteq L\backslash
\{1\}$. Then $\mathcal{L}$ is an $X$-top lattice if and only if every $x\in
X $ is strongly $C^{X}(L)$-irreducible in $(C^{X}(L),\wedge )$ (i.e. $%
X=SI^{C^{X}(L)}(X)$).
\end{thm}

The following is a direct, but very useful, consequence of Theorem \ref{xct}
especially in constructing examples and counterexamples.

\begin{cor}
\label{Y-top}Let $\mathcal{L}=(L,\wedge ,\vee ,1,0)$ be an $X$-top lattice
for some $X\subseteq L\backslash \{1\}$. If $Y\subseteq X,$ then $L$ is a $Y$%
-top lattice and the corresponding topology on $Y$ is the induced subspace
topology.
\end{cor}

\begin{Beweis}
We show first that $C^{Y}(L)\subseteq C^{X}(L)$. Let $a\in C^{Y}(L).$ Notice
that $V_{Y}(a)\subseteq V_{X}(a)$ as $Y\subseteq X,$ hence, $\bigwedge
V_{X}(a)\leq \bigwedge V_{Y}(a)=a$. Since, $a\leq \bigwedge V_{X}(a)$, we
conclude that $a=\bigwedge V_{X}(a)$. By the direct implication of Theorem %
\ref{xct}, $Y\subseteq SI^{C^{X}(L)}(X).$ Since $C^{Y}(L)\subseteq C^{X}(L),$
and both share the \emph{same meet}, we conclude that $Y\subseteq
SI^{C^{Y}(L)}(Y).$ By the reverse implication of Theorem \ref{xct}, $L$ is a
$Y$-top lattice. Notice that $V_{Y}(a)=V_{X}(a)\cap Y$ for any $a\in L$,
i.e. the corresponding topology on $Y$ is the induced subspace topology on $%
Y\hookrightarrow X.\blacksquare $
\end{Beweis}

\begin{remarks}
\label{finite-BMax}Let $\mathcal{L}=(L,\wedge ,\vee ,1,0)$ be an $X$-top
lattice for some $X\subseteq L\backslash \{1\}$.

\begin{enumerate}
\item $X=SI^{X}(X).$

\item If $X$ is finite, then $X=CSI^{X}(X).$

\item If $Max(X)\ $is finite, then $X$ is BMax.

\item If $Min(X)\ $is finite, then $X$ is AMin.
\end{enumerate}
\end{remarks}

\subsection*{Semirings}

We recall some definitions and examples from the Theory of Semirings. We
follow \cite{Gol1999} (unless otherwise stated explicitly).

\begin{defn}
\cite{Gol1999} A \textbf{semiring} $R$ is a non-empty set equipped with two
binary operations: addition "+" and multiplication "$\cdot $" that satisfy
the following conditions:

\begin{enumerate}
\item $(R,+)$ is a commutative monoid with neutral element $0$;

\item $(R,\cdot )$ is a monoid with neutral element $1$;

\item Multiplication distributes over addition from either side;

\item $0\cdot r=0=r\cdot 0$ for all $r\in R$ (i.e. $0$ is \emph{absorbing}).

\item $1\neq 0$.

If, in addition, $(R,\cdot )$ is commutative, we say that $R$ is a \emph{%
commutative semiring}. We call a commutative semiring with no non-zero
zero-divisors a \textbf{semidomain}. A semiring that is not a ring is called
a \textbf{proper semiring.}
\end{enumerate}
\end{defn}

\begin{punto}
Let $(S,+,0,\cdot ,1)$ be a semiring. We call $a\in S$ \textbf{additively
idempotent} (resp. \textbf{multiplicatively idempotent}) iff $a+a=a$ (resp. $%
a^{2}=a$). We call $a\in S$ \textbf{idempotent} iff $a$ is additively
idempotent and multiplicatively idempotent. With\ $I^{+}(S)$ (resp. $%
I^{\times }(S),$ $I(S)$) we denote the subset of additively idempotent
(resp. multiplicatively idempotent, idempotent) elements of $S.$ We say that
$S$ is an \textbf{additively idempotent semiring} (resp. \textbf{%
multiplicatively idempotent} \textbf{semiring}, \textbf{idempotent semiring}%
) iff $I^{+}(S)=S$ (resp. $I^{\times }(S)=S,$ $I(S)=S$).
\end{punto}

\begin{punto}
Let $R$ be a semiring. We say that $A\subseteq R$ is \textbf{subtractive}
iff for all $r\in R$ and $a\in A,$ we have: $r+a\in A\Longrightarrow r\in A.$
We call a left (right) ideal of $R$ a left (right) $k$\textbf{-ideal}, iff $%
I\subseteq R$ is subtractive. We say that $R$ is a \textbf{left subtractive
semiring} (\textbf{right subtractive semiring}) iff every left (right) ideal
of $R$ is subtractive. A \textbf{subtractive semiring} is a semiring that is
left subtractive and right subtractive.
\end{punto}

\begin{exs}
\begin{enumerate}
\item Every ring is a semiring.

\item $(\mathbb{W},+,0,\cdot ,1)$, where $\mathbb{W}=\{0,1,2,3,\cdots \}$ is
the set of \emph{whole numbers} (non-negative integers) is a commutative
semiring.

\item Every \emph{bounded distributive lattice} is an idempotent commutative
semiring.

\item The \textbf{Boolean semiring }is $\mathbb{B}=\{0,1\}$, in which $%
1+1=1. $
\end{enumerate}
\end{exs}

\begin{punto}
Let $R$ be a (semi)ring and consider $Ideal(R),$ the complete lattice of all
ideals of $R$ with $I\vee J:=I+J$ and $I\wedge J:=I\cap J$ for ideals $I,J$
of $R.$ If $X\subseteq Ideal(R)\backslash \{R\},$ then we say that $R$ is an
$X$\textbf{-top (semi)ring} iff $Ideal(R)$ is an $X$-top lattice. For
example, $R$ is a $Spec(R)$-top (semi)ring and the topology on $Spec(R)$ is
the ordinary \emph{Zariski topology} on the spectrum of prime ideals of $R$
(e.g. \cite{AM1969}). Moreover, for any subset $Y\subseteq Spec(R)$, we
have: $R$ is a $Y$-top (semi)ring by Corollary \ref{Y-top}. In particular, $%
R $ is a $Max(R)$-top (semi)ring and a $Min(R)$-top (semi)ring.$\blacksquare
$
\end{punto}

\subsection*{General Topology}

In what follows, we recall some definitions and elementary results from
General Topology. We follow \cite{Wil1970} unless otherwise mentioned
explicitly.

\bigskip

\begin{notatation}
Let $X$ be a topological space. For $Y\subseteq X,$ we denote with $\mathcal{%
O}(Y)$ (resp. $\mathcal{C}(Y),$ $\mathcal{CO}(Y)$, $\mathcal{CD}(Y),$ $%
\mathcal{K}(Y),$ $\mathcal{KO}(Y),$ $\mathcal{KC}(Y),$ $\mathcal{KCO}(Y),$ $%
\mathcal{I}(Y),$ $\mathcal{IC}(Y)$) the collection of all open (resp.
closed, clopen, connected, compact, compact open, compact closed, compact
clopen, irreducible closed, maximal irreducible closed) subsets of $X$ that
contain $Y.$ By abuse of notation (compare with the notation for $%
Y=\emptyset $), we denote with $\mathcal{O}(X)$ (resp. $\mathcal{C}(X),$ $%
\mathcal{CO}(X)$, $\mathcal{CD}(X),$ $\mathcal{K}(X),$ $\mathcal{KO}(X),$ $%
\mathcal{KC}(X),$ $\mathcal{KCO}(X),$ $\mathcal{I}(X),$ $\mathcal{IC}(X)$)
the collection of all open (resp. closed, clopen, connected, compact,
compact open, compact closed, compact clopen, irreducible closed, maximal
irreducible closed) subsets of $X.$
\end{notatation}

\subsubsection*{Separation Axioms}

\begin{defns}
We say that a topological space $X$ is

\begin{enumerate}
\item $T_{0}$ (\textbf{Kolmogorov)} iff any two distinct points $x\neq y$ in
$X$ are \emph{topologically distinguishable,} i.e. iff there exists an open
neighborhoods $U_{x}$ of $x$ such that $y\notin U_{x}$ \emph{or} an open
neighborhood of $U_{y}$ of $y$ such that $x\notin U_{y}.$

\item $R_{0}$ (\textbf{symmetric})\ iff any two topologically
distinguishable points $x\neq y$ in $X$ are \emph{separated} (i.e. there
exist two \emph{open }neighborhoods $U_{x}$ of $x$ and $U_{y}$ of $y$, such
that $y\notin U_{x}$ and $x\notin U_{y}$);

\item $T_{1}$ (\textbf{Fr\'{e}chet)} iff any $x\neq y$ in $X$ are \emph{%
separated }(\emph{i.e. }$T_{1}=T_{0}+R_{0}$);

\item $R_{1}$ (\textbf{preregular}) iff any two topologically
distinguishable points of $X$ can be separated by \emph{disjoint open}
neighborhoods;

\item $T_{1\frac{1}{2}}$ (\textbf{KC}) iff all compact subsets of $X$ are
closed;

\item $T_{2}$ (\textbf{Hausdorff}) iff any $x\neq y$ in $X$ are \emph{%
separated by disjoint open neighborhoods }($T_{2}=T_{0}+R_{1}=T_{1}+R_{1}$);

\item \textbf{connected} iff $X$ cannot be written as the union of two \emph{%
disjoint} proper closed (open) subsets;

\item \textbf{irreducible} (\textbf{hyperconnected}) iff $X$ cannot be
written as the union of two \emph{proper} closed subsets (\emph{whether
disjoint or non-disjoint}).
\end{enumerate}
\end{defns}

\begin{rem}
\label{anti-hyper}We say that a topological space is $X$ is \textbf{anti-}$%
T_{2}$ (\textbf{anti-Hausdorff}) iff $\left\vert X\right\vert \geq 2$ and no
two distinct points of $X$ can be separated by disjoint open sets. By \cite[%
Theorem 4.2]{MM2012}, $X$ is anti-$T_{2}$ if and only if $X$ is irreducible
(hyperconnected)\textbf{.}
\end{rem}

\begin{defn}
Let $X$ be a topological space. For each $x\in X,$ the \textbf{connected
component} of $C(x)$ of $x$ and the \textbf{quasicomponent} $Q(x)\ $of $x$
are%
\begin{equation*}
C(x)=\bigcup\limits_{C\in \mathcal{CD}(x)}C\text{ and }Q(x):=\bigcap_{W\in
\mathcal{CO}(x)}W.
\end{equation*}%
We say that $X$ is \textbf{totally disconnected} (resp. \textbf{totally
separated) }iff $C(x)=\{x\}$ (resp. $Q(x)=\{x\}$) for every $x\in X.$
\end{defn}

\begin{remarks}
\label{CO-sep}Let $X$ be a topological space.

\begin{enumerate}
\item $X$ is totally separated if and only if whenever $x\neq y$ in $X,$
there is $U\underset{\text{clopen}}{\subseteq }X$ such that $x\in U$ while $%
y\notin U$ (and $V\underset{\text{clopen}}{\subseteq }X$ such that $y\in V$
while $x\notin V$). In particular, every totally separated topological space
is $T_{2}.$

\item If $X$ is a totally disconnected topological space, then $X$ is $%
T_{1}. $ Let $x\in X$. Since $\{x\}$ is connected, and \emph{the closure of
any connected subset of a topological space is connected}, it follows that $%
\overline{\{x\}}\subseteq C(x)=\{x\}$, whence $\overline{\{x\}}=\{x\}.$
\end{enumerate}
\end{remarks}

\section{$X$-Top lattices with the $T_{1}$ Separation Axiom}

\qquad In this section, we study the $X$-top lattices for which $X$ is $%
T_{1}.$

\begin{punto}
Let $L$ be an $X$-top lattice for some $X\subseteq L\backslash \{1\}.$ Then $%
(X,\leq )$ is a poset. We define the \textbf{height of} $x\in X$ as%
\begin{equation*}
ht(x):=\sup \{n\geq 0\text{ }\mid \text{there exists }\{x_{0},\cdots
,x_{n}\}\subseteq X\text{ with }x_{0}\lneqq \cdots \lneqq x_{n}=x\}
\end{equation*}%
and the \textbf{Krull dimension of }$X$ as $K.\dim (X)=\sup \{ht(x)\mid x\in
X\}.$
\end{punto}

\begin{punto}
Let $R$ be a commutative (semi)ring, $L=Ideal(R)$ and $X=Spec(R).$ Then the
\emph{Krull dimension} $K.\dim (R)$ of $R$ is nothing but $K.\dim (Spec(R)).$
For example, we have $K.\dim (\mathbb{Z})=K.\dim (Spec(\mathbb{Z}))=1$ and $%
K.\dim (Max(\mathbb{Z}))=0=K.\dim (Max(\mathbb{Z})).$
\end{punto}

\begin{defn}
(cf. \cite{G-R2007}, \cite{Dun1977}, \cite{Wil1970}) Let $X$ be a
topological space and set for $x\in X:$%
\begin{equation*}
Ker(x):=\bigcap\limits_{U\in \mathcal{O}(x)}U\text{ and }E(x):=\bigwedge%
\limits_{y\in X\backslash \{x\}}y.
\end{equation*}%
We say that $x\in X$ is

\begin{enumerate}
\item \textbf{isolated} iff $\{x\}$ is an open set;

\item \textbf{kerneled} iff $\{x\}=Ker(x);$

\item \textbf{regular open} iff $\{x\}=int(\overline{\{x\}});$

\item \textbf{excluded} iff $E(x)=\bigwedge D(x).$
\end{enumerate}
\end{defn}

\bigskip

For any topological space $X,$ we set%
\begin{equation*}
\begin{array}{lllllll}
K(X) & := & \{x\in X\text{ }\mid \text{ }\{x\}=Ker(\{x\})\}; &  & Iso(X) & :=
& \{x\in X\text{ }|\text{ }\{x\}\text{ is open}\}; \\
RO(X) & := & \{x\in X\text{ }|\text{ }\{x\}=int(\overline{\{x\}})\}; &  &
Cl(X) & := & \{x\in X\text{ }|\text{ }\{x\}\text{ is closed}\}; \\
\mathcal{E}(X) & := & \{x\in X\mid E(x)=\bigwedge D(x)\}. &  &  &  &
\end{array}%
\end{equation*}

\bigskip

\begin{lem}
\label{kerneled}Let $L$ be an $X$-top lattice for some $X\subseteq
L\backslash \{1\}.$

\begin{enumerate}
\item $Cl(X)=Max(X);$

\item $K(X)=Min(X);$

\item $RO(X)\subseteq Iso(X)\subseteq Min(X).$
\end{enumerate}
\end{lem}

\begin{Beweis}
\begin{enumerate}
\item This follows from the fact that for any $x\in X,$ we have%
\begin{equation*}
\overline{\{x\}}=V(x):=\{y\in X\mid x\leq y\}.
\end{equation*}

\item This follow from the fact that for any $x\in X,$ we have%
\begin{equation*}
\begin{tabular}{lllll}
$Ker(\{x\})$ & $=$ & $\bigcap\limits_{U\in \mathcal{O}(x)}U$ & $=$ & $%
\bigcap \{D(a)\ |\ x\in D(a)\}$ \\
& $=$ & ${\bigcap \{X\backslash V(a)\ |\ a>x\ \text{or}\ a||x\}}$ & $=$ & $%
\{z\in X\ |\ z\leq x\}$%
\end{tabular}%
\end{equation*}

\item This follows from the definitions and (2).$\blacksquare $
\end{enumerate}
\end{Beweis}

\begin{defn}
A topological space $X$ is said \textbf{sober} iff every irreducible closed
subset $Y\subseteq X$ has a \emph{unique} \textbf{generic point} (i.e. $y\in
Y$ such that $Y=\overline{\{y\}}$).
\end{defn}

\begin{defn}
(cf. \cite{Est1988}, \cite{Hoc1969}, \cite{DST2019}) A \textbf{spectral space%
} is a topological space $X$ that satisfies any (hence all) of the following
equivalent conditions:

\begin{enumerate}
\item $X$ is sober, compact and has a base $\mathcal{B}\subseteq \mathcal{KO}%
(X)$ closed under finite intersections;

\item $X$ is homeomorphic to $Spec(R)$ for some \emph{commutative (semi)ring}
$R;$

\item $X$ is homeomorphic to a projective limit of \emph{finite }$T_{0}$
spaces.
\end{enumerate}
\end{defn}

\begin{defn}
A topological space is \textbf{quasi-Hausdorff} \cite{Hoc1969} iff for any $%
x\neq y$ in $X:$ \emph{either} $x$ and $y$ are separated by disjoint open
neighborhoods\emph{\ or }there exists $z\in X$ such that $\{x,y\}\subseteq
\overline{\{z\}}.$
\end{defn}

\begin{lem}
\label{T1-qH-T2}

\begin{enumerate}
\item Every finite $T_{0}$ space is spectral \emph{(cf. \cite{Hoc1969})}.

\item Every spectral space is quasi-Hausdorff \emph{(cf. \cite[Corollary 2,
page 45]{Hoc1969})}.

\item A topological space $X$ is $T_{2}$ if and only if $X$ is $T_{1}$ and
quasi-Hausdorff.

\item A $T_{1}$ spectral space is $T_{2}$ (\emph{cf. \cite[Exercise 3.11]%
{AM1969}}).
\end{enumerate}
\end{lem}

\begin{prop}
\label{dim-0}Let $\mathcal{L}=(L,\wedge ,\vee ,1,0)$ be an $X$-top lattice
for some $X\subseteq L\backslash \{1\}.$

\begin{enumerate}
\item $X$ is $T_{0}.$

\item If $X$ is finite, then $X$ is spectral.

\item $Max(X)=Max(C^{X}(L)).$

\item $X$ is $T_{1}\Longleftrightarrow X$ is $R_{0}\Longleftrightarrow $ $%
K.\dim (X)=0.$

\item $X$ is $T_{2}\Longleftrightarrow X$ is $R_{1}\Longleftrightarrow
K.\dim (X)=0$ and $X$ is quasi-Hausdorff.
\end{enumerate}
\end{prop}

\begin{Beweis}
\begin{enumerate}
\item Let $x\neq y$ in $X.$ Then either $x\nleqq y$ whence $y\in
D(x):=X\backslash V(x)$ \emph{or} $y\nleqq x$ whence $x\in D(y):=X\backslash
V(y)$.

\item Every finite $T_{0}$ space is spectral (cf. \cite{Hoc1969}).

\item Notice that $X\subseteq C^{X}(L)\backslash \{1\}.$

Suppose that $m\in Max(X)\backslash Max(C^{X}(L)).$ Then there exists some $%
q\in C^{X}(L)$ such that $m\nleqq q\nleqq 1.$ Since $V(q)\neq \{1\},$ there
exists some $x\in X\subseteq L\backslash \{1\}$ with $q\leq x,$ whence $%
m\nleqq x\nleqq 1$ (a contradiction).

On the other hand, let $y\in Max(C^{X}(L)).$ Then $V(y)=\{y\}$ (otherwise, $%
y=\bigwedge\limits_{x\in V(y)}<z$ for some $z\in X\subseteq C^{X}(L)$ with $%
y<z<1,$ a contradiction). So, $y\in Max(X).$

\item The first equivalent follows from (1) and the fact that $%
T_{1}=T_{0}+R_{0}.$ The second equivalence follows from%
\begin{equation*}
X\text{ is }T_{1}\Longleftrightarrow X=Cl(X)\overset{\text{Lemma \ref%
{kerneled}}}{\Longleftrightarrow }X=Max(X)\Longleftrightarrow
Min(X)=X=Max(X).
\end{equation*}

\item The first equivalence follows from (1) and the fact that $%
T_{2}=T_{0}+R_{1}.$ The second equivalence follows by (4)\ and Lemma \ref%
{T1-qH-T2} (3).$\blacksquare $
\end{enumerate}
\end{Beweis}

\begin{defn}
We call a (semi)ring $R$

\textbf{von Neumann regular} iff for every $a\in R$ there exists $b$ such
that $a=aba;$

$\pi $\textbf{-regular} iff for every $a\in R$ there exists $b\in R$ and $%
n\geq 1$ such that $a^{n}=a^{n}ba^{n};$

\textbf{semiprime }iff $P(R):=\bigcap\limits_{P\in Spec(R)}P=0;$

\textbf{reduced }iff $Nil(R):=\{a\in R$ $\mid $ $a^{n}=0$ for some $n\geq
1\}=0.$

If $R$ is commutative, then $P(R)=Nil(R)$ (cf. \cite[Proposition 1.8]{AM1969}%
), whence a semiprime commutative (semi)ring is the same as a reduced one.
\end{defn}

\begin{defn}
(cf. \cite[29.4]{Wil1970})\ We say that a topological space $X$ is \textbf{%
inductively zero-dimensional }iff $X$ has a base of \emph{clopen sets}
(equivalently, $ind.\dim (X)=0,$ where $ind.\dim (X)$ is the so-called
\textbf{small inductive dimension} of $X$ \cite[page 105]{AP1990}).
\end{defn}

\begin{rem}
\label{z-cr}Every inductively zero-dimensional topological space is
completely regular (e.g. \cite[Proposition I.31]{AP1990}). So, every $T_{0}$
inductively zero-dimensional topological space is $T_{3\frac{1}{2}}.$
\end{rem}

\begin{lem}
\label{T1-td}Let $X$ be an \emph{inductively zero-dimensional }topological
space. The following are equivalent:

\begin{enumerate}
\item $X$ is totally separated;

\item $X$ is totally disconnected;

\item $X$ is $T_{1};$

\item $X$ is $T_{0};$

\item $X$ is $T_{3\frac{1}{2}};$

\item $X$ is $T_{2}.$
\end{enumerate}
\end{lem}

\begin{Beweis}
The implications $(5\Longrightarrow 6\Longrightarrow 3\Longrightarrow 4)$\
are trivial.

$(1\Longrightarrow 2)$ Every totally separated space is totally disconnected
(since $C(x)\subseteq Q(x)$ for every point $x$ in the space).

$(2\Longrightarrow 3)$ Every totally disconnected space is $T_{1}$ by Remark %
\ref{CO-sep}.

$(3\Longrightarrow 1)$ Let $X$ be inductively zero-dimensional with a base $%
\mathcal{B}\subseteq \mathcal{CO}(X)$ and assume that $X$ is $T_{1}.$ Let $%
x\in X$ and pick any $y\neq x$ in $X.$ Since $X$ is $T_{1},$ there exists
some $U\in \mathcal{O}(x)$ such that $x\in U$ and $y\notin U.$ Since $%
\mathcal{B}$ is a base, there exists $\widetilde{U}\in \mathcal{B}$ such
that $x\in \widetilde{U}\subseteq U.$ Since $\widetilde{U}$ is clopen, we
have $\widetilde{U}\subseteq C(x)$ (otherwise $C(x)=(C(x)\cap \widetilde{U}%
)\cup (C(x)\cap X\backslash \widetilde{U})$ a union of two \emph{disjoint}
proper closed subsets of $C(x)$, a contradiction). So, $y\notin Q(x).$ We
conclude that $Q(x)=\{x\}.$ Since $x\in X$ was arbitrary, we conclude that $%
X $ is totally separated.

$(4\Longrightarrow 5)$ This follows by Remark \ref{z-cr}.$\blacksquare $
\end{Beweis}

\begin{defn}
(cf. \cite{Joh1982}) A \textbf{Stone space (Boolean space, profinite space})
is a topological space $X$ that satisfies any, hence all, of the following
equivalent conditions:

\begin{enumerate}
\item $X$ is $T_{0},$ inductively zero-dimensional space and compact;

\item $X$ is homeomorphic to a projective limit of \emph{finite discrete}
spaces.
\end{enumerate}
\end{defn}

\begin{prop}
\label{abs-flat}Let $L$ be an $X$-top space for some $X\subseteq L\backslash
\{1\}.$ The following are equivalent:

\begin{enumerate}
\item $X$ is a Stone space;

\item $X$ is spectral and inductively zero-dimensional;

\item $X$ is spectral and totally separated;

\item $X$ is spectral and $T_{2};$

\item $X$ is spectral and $T_{1\frac{1}{2}};$

\item $X$ is spectral and $T_{1};$

\item $X$ is spectral and $K.\dim (X)=0;$

\item $X$ is homeomorphic to $Spec(R)$ for some Krull $0$-dimensional
commutative (semi)ring $R;$

\item $X$ is homeomorphic to $Spec(R)$ for some commutative ring $R$ with $%
R/Nil(R)$ von Neumann regular;

\item $X$ is homeomorphic to $Spec(R)$ for some $\pi $-regular commutative
ring $R.$
\end{enumerate}
\end{prop}

\begin{Beweis}
$(1\Longrightarrow 2)$ Every Stone space is spectral (being the projective
limit of finite discrete, whence finite $T_{0},$ topological spaces).

$(2\Longrightarrow 3)$ This follows by Lemma \ref{T1-td}.

$(3\Longrightarrow 4)\ $Every totally separated topological space is $T_{2}$
by Lemma \ref{CO-sep}.

$(4\Longrightarrow 5)\ $Every $T_{2}$ space is $T_{1\frac{1}{2}}$ by \cite[%
Theorem 17.5 (2)]{Wil1970}.

$(5\Longrightarrow 1)$ Since $X$ is spectral, $\mathcal{KO}(X)$ is a base of
$X.$ Since $X$ is $T_{1\frac{1}{2}},$ every compact subset of $X$ is closed,
whence $\mathcal{KO}(X)$ is a base of (compact) clopen sets. So, $X$ is $%
T_{0},$ compact and inductively zero-dimensional.

$(4\Longleftrightarrow 6)$ This follows by Lemma \ref{T1-qH-T2} (4).

$(6\Longleftrightarrow 7)$ This follows by Proposition \ref{dim-0} (4).

$(7\Longleftrightarrow 8)$ This follows by the characterization of a
spectral space as one homeomorphic to the prime spectrum of a commutative
(semi)ring (\cite{Est1988}, \cite{Hoc1969}). The Krull dimension is clearly
preserved by any homeomorphism.

$(8\Longleftrightarrow 9\Longleftrightarrow 10)$ These follow by well-known
characterizations of Krull $0$-dimensional commutative rings (e.g. \cite%
{Sto1968}, \cite[Exercise 3.11]{AM1969}, \cite{Gil2000}).$\blacksquare $
\end{Beweis}

The following result is a generalization of \cite[Theorem 3.1]{PV2011} (and
\cite[Theorem 5.1]{PRV2009}):

\begin{prop}
\label{iso}Let $L$ be an $X$-top lattice for some $X\subseteq L\backslash
\{1\}.$ Then%
\begin{equation*}
Iso(X)=Min(X)\cap CSI^{X}(X).
\end{equation*}
\end{prop}

\begin{Beweis}
$(\implies )$ Notice that $Iso(X)\subseteq K(X)=Min(X)$ by Lemma \ref%
{kerneled} (2). Let $x\in Iso(X)$ and suppose that $x\notin CSI^{X}(X).$
Then there exists $A\subseteq X$ such that $\bigwedge\limits_{a\in A}a\leq x$
but $a\nleq x$ for all $a\in A.$ Since $x$ is isolated, $\{x\}=D(z)=X%
\backslash V(z)$ for some $z\in L.$ Notice that $A\subseteq V(z)$. But then,
$z\leq \bigwedge\limits_{a\in A}a\leq x,$ whence $x\in V(z)$, a
contradiction. Therefore, $x\in CSI^{X}(X).$

$(\impliedby )$ Let $x\in Min(X)\cap CSI^{X}(X).$ Set ${z:=\bigwedge%
\limits_{y\in X\backslash \{x\}}y}$. Since $x$ is minimal and completely
strongly $X$-irreducible, $z\nleq x,$ i.e. $x\in D(z).$ Observe that, by our
choice of $z,$ we have $X\backslash \{x\}\subseteq V(z).$ So, $D(z)=\{x\},$
i.e. $x\in Iso(X).\blacksquare $
\end{Beweis}

The following result recovers \cite[Proposition 1.23, Proposition 1.24]%
{AH2019} and gives additional characterizations of an $X$-top lattice for
which $X$ is discrete.

\begin{thm}
\label{discrete}Let $L$ be an $X$-top lattice for some $X\subseteq
L\backslash \{1\}$. The following are equivalent:

\begin{enumerate}
\item $X$ is discrete;

\item Every $x\in X$ is minimal and completely strongly $X$-irreducible;

\item $X$ is PAMin (every $x\in X$ is absolutely $X$-minimal);

\item Every $x\in X$ is maximal and $X$ is BMax;

\item $X$ is PBMax (every $x\in X$ is barely $X$-maximal);

\item $X$ is $T_{1}$ and $X$ is BMax;

\item $X$ is $T_{1}$ and $C^{X}(L)$ has the complete max-property.

In this case, we have%
\begin{equation*}
AMin(X)=Min(X)=X=Max(X)=BMax(X).
\end{equation*}
\end{enumerate}
\end{thm}

\begin{Beweis}
Notice that under any of these conditions, we can assume $X$ to be $T_{1}$
(by Proposition \ref{dim-0}) and that%
\begin{equation*}
Min(X)=X=Max(X)=Max(C^{X}(L)).
\end{equation*}

$(1)\Longleftrightarrow (2):$ This follows since%
\begin{equation*}
X\text{ is discrete }\Longleftrightarrow X=Iso(X)\overset{\text{Proposition %
\ref{iso}}}{=}Min(X)\cap CSI^{X}(X).
\end{equation*}

$(2)\Longleftrightarrow (3)$ We have $X=CSI^{Min(X)}(Min(X))=AMin(X).$

$(2)\Longleftrightarrow (4)\ $This follows since $Min(X)=X=Max(X).$

$(4)\Longleftrightarrow (5)\ $We have $%
CSI^{X}(X)=CSI^{Max(X)}(Max(X))=BMax(X).$

$(5)\Longleftrightarrow (6)$ Notice that $X=Max(X).$ So, $X$ has the \emph{%
complete max property} if and only if $X=BMax(X).$

$(6)\Longleftrightarrow (7)$ This follows since $Max(X)=Max(C^{X}(L))$ by
Proposition \ref{dim-0} (3).$\blacksquare $
\end{Beweis}

\begin{punto}
We call a bounded join-semilattice $(L,\vee ,1)$ \textbf{finitely joinable},
iff whenever $1=\bigvee\limits_{\lambda \in \Lambda }a_{\lambda }$ for $%
\{a_{\lambda }\mid \lambda \in \Lambda \}\subseteq L,$ there exists a \emph{%
finite} subset $\{\lambda _{1},\cdots ,\lambda _{n}\}\subseteq \Lambda $
such that $1=\bigvee\limits_{i=1}^{n}a_{\lambda _{i}}.$ For example, $%
L=(Ideal(R),+,R)\ $is finitely joinable for every (semi)ring $R$ with $%
1_{R}. $
\end{punto}

\begin{cor}
\label{X-Max-discrete}Let $L$ be an $X$-top lattice for some $X\subseteq
L\backslash \{1\}$. Then $L$ is $Max(X)$-top, $Max(X)$ is $T_{1}$ and the
following are equivalent:

\begin{enumerate}
\item $X$ is BMax;

\item $J(X):=\bigwedge\limits_{\mathfrak{m}\in Max(X)}\mathfrak{m}$ is an
irredundant intersection of maximal elements of $X;$

\item $Max(X)$ is discrete;

If $L\backslash \{1\}$ is $X$-coatomic and $L$ is finitely joinable, then
these are equivalent to:

\item $X$ is FMax.
\end{enumerate}
\end{cor}

\begin{Beweis}
Notice that $L$ is $Max(X)$-top by Corollary \ref{Y-top}. It's obvious that $%
\dim (Max(X))=0$ (equivalently, $Max(X)\ $is $T_{1}$ by Proposition \ref%
{dim-0}).

$(1\Longleftrightarrow 2)$ is obvious.

$(1\Longleftrightarrow 3)$ Apply Theorem \ref{discrete} (or Proposition \ref%
{iso}) to $Max(X)\subseteq L\backslash \{1\}.$

$(4\Longrightarrow 1)$ Since $L$ is an $X$-top lattice, this follows by
Remark \ref{finite-BMax} (3).

Let $L\backslash \{1\}\ $be $X$-coatomic and $L$ be finitely joinable.
\textbf{Claim:} $(3\Longrightarrow 4).$

Assume that $Max(R)$ is discrete. For each $\mathfrak{m}\in Max(X),$ pick
some $a_{\mathfrak{m}}\in L$ such that $D(a_{\mathfrak{m}})\cap Max(X)=\{%
\mathfrak{m}\}.$ Consider $a:=\bigvee\limits_{\mathfrak{m}\in Max(R)}a_{%
\mathfrak{m}}.$ If $a\neq 1,$ then (since $L\backslash \{1\}$ is $X$%
-coatomic) there exists some $\mathfrak{m}_{a}\in Max(X)$ such that $a\leq
\mathfrak{m}_{a}$ whence $a_{\mathfrak{m}_{a}}\leq \mathfrak{m}_{a},$ a
contradiction. So, $a=1.$ Since $1$ is \emph{finitely joinable}, there
exists $\{\mathfrak{m}_{1},\cdots ,\mathfrak{m}_{n}\}\subseteq Max(X)\ $such
that $1=\bigvee\limits_{i=1}^{n}$ $a_{\mathfrak{m}_{i}}.$ It follows that $%
Max(X)\subseteq \bigcup\limits_{i=1}^{n}(D(a_{\mathfrak{m}_{i}})\cap
Max(X))=\{\mathfrak{m}_{1},\cdots ,\mathfrak{m}_{n}\}.$ So, $\left\vert
Max(X)\right\vert <\infty ,$ i.e. $X$ is FMax.$\blacksquare $
\end{Beweis}

\begin{cor}
\label{X-Min-discrete}Let $L$ be an $X$-top lattice for some $X\subseteq
L\backslash \{1\}.$ Then $L$ is $Min(X)$-top and $Min(X)$ is $T_{1}$ and the
following are equivalent:

\begin{enumerate}
\item $X$ is AMin;

\item $Min(X)$ is discrete;

\item $Q(X):=\bigwedge\limits_{q\in Min(X)}q$ is an irredundant meet.

If $X$ is atomic and $P(X):=\bigwedge\limits_{x\in X}q,$ then these are
equivalent to

\item $P(X)=\bigwedge\limits_{q\in Min(X)}q$ is an irredundant intersection.
\end{enumerate}
\end{cor}

\begin{rem}
\label{Spec-discrete}Every spectral space is $T_{0}$ (since it is sober) and
compact. On the other hand, every finite topological space is spectral (e.g.
\cite{Hoc1969}). If $L$ is an $X$-top lattice for some $X\subseteq
L\backslash \{1\},$ then the following are equivalent by Proposition \ref%
{dim-0}:

\begin{enumerate}
\item $X$ is spectral and discrete;

\item $X$ is compact and discrete;

\item $X$ is finite and discrete;

\item $X$ is finite and $T_{1}$ ($R_{0},$ $R_{1},$ $T_{2}$);

\item $Max(X)$ is finite and $X$ is $T_{1};$

\item $Max(X)$ is finite and $K.\dim (X)=0;$

\item $Min(X)$ is finite and $X$ is $T_{1};$

\item $Min(X)$ is finite and $K.\dim (X)=0.$
\end{enumerate}
\end{rem}

\begin{defn}
We say that a (semi)ring $R$ is

\textbf{FMax-(semi)ring} iff $R$ has finitely many maximal ideals;

\textbf{BMax-(semi)ring }iff every maximal ideal of $R$ is barely maximal;

\textbf{PBMax-(semi)ring }iff every prime ideal of $R$ is barely maximal;

\textbf{AMin-(semi)ring }iff every minimal prime ideal of $R$ is absolutely
minimal;

\textbf{FMin-(semi)ring} iff $R$ has finitely many minimal prime ideals;

\textbf{PAMin-(semi)ring }iff every prime ideal of $R$ is absolutely minimal.
\end{defn}

\begin{rem}
A (semi)ring $R$ is a PAMin-(semi)ring $\Longleftrightarrow $ $K.\dim (R)=0$
and $R$ is a AMin-(semi)ring $\Longleftrightarrow $ $K.\dim (R)=0$ and $R$
is a BMax-(semi)ring $\Longleftrightarrow $ $R$ is a PBMax-(semi)ring.
\end{rem}

Let $R$ be a commutative (semi)ring and consider the complete lattice $%
(Ideal(R),\cap ,+,R,0).$ By \cite[Proposition 6.59]{Gol1999}, every ideal of
$R$ is contained in a maximal ideal of $R$ which is prime by \cite[Crollary
7.13]{Gol1999}, i.e. $Ideal(R)\backslash \{R\}$ is $Spec(R)$-coatomic.\
Clearly, $R$ is finitely joinable as $1_{R}\in R.$ Taking these observations
into account, we obtain as a direct consequence of Corollary \ref%
{X-Max-discrete}:

\begin{cor}
\label{semilocal-discrete}The following are equivalent for a commutative
(semi)ring $R:$

\begin{enumerate}
\item $R$ is a BMax-(semi)ring;

\item $J(R)$ is an irredundant intersection of maximal ideals of $R;$

\item $Max(R)$ is discrete;

\item $R$ is an FMax (semi)ring.
\end{enumerate}
\end{cor}

\begin{defn}
(\cite{Wis1991}, \cite{Lam2001}) We say that a (semi)ring $R$ is \textbf{%
local} iff $R$ has a unique maximal left ideal.
\end{defn}

\begin{rem}
A ring $R$ is said to be \textbf{semilocal} iff $R/Jac(R)$ is semisimple
(e.g. \cite{Wis1991}). An FMax-ring is semilocal and the converse holds for
\emph{commutative} rings (e.g. \cite[Proposition 20.2]{Lam2001}). Let $R$ be
a commutative ring. We say that $R$ has the (\emph{complete})\ \emph{%
max-property} iff the lattice $(Ideal(R),\cap ,+,R,0)$ has the (\textbf{%
complete}) \textbf{max-property. }Notice that for $X=Spec(R),$ we have $%
Max(L)=Max(X).$ Every ring $R$ has the max-property. A commutative ring $R$
has the \emph{complete max-property} (equivalently, $X$ is BMax) if and only
if $R$ is a \emph{semilocal ring} as one can notice by combining \cite[%
Theorems 4.4 and 6.8]{Smi2011}. We recovered this result and extended it to
arbitrary commutative semirings in Corollary \ref{semilocal-discrete} as a
consequence of Corollary \ref{finite-BMax} which has a very short proof.
\end{rem}

\begin{ex}
\label{Zn}Let $n\geq 2$ with prime factorization $n=p_{1}^{m_{1}}\cdots
p_{k}^{m_{k}}.$ Then $\mathbb{Z}_{n}$ is semilocal and $Spec(\mathbb{Z}%
_{n}))=Max(\mathbb{Z}_{n})=\{(p_{1}),\cdots ,(p_{n})\}.$ One can check
easily that $\mathbb{Z}_{n}$ is a BMax-ring: $\bigcap\limits_{j=1,i\neq
j}^{k}(p_{j})\nsubseteqq (p_{i})$ for any $i=1,\cdots ,k$ and that $Spec(%
\mathbb{Z}_{n})$ is discrete: $\{(p_{i})\}=D((\frac{n}{p_{i}^{m_{i}}}))$ for
$i=1,\cdots ,k.\blacksquare $
\end{ex}

By \cite[Proposition 7.14]{Gol1999}, every prime ideal $P$ of a commutative
semiring $R$ contains a minimal prime ideal $Q$ (i.e. $Spec(R)$ is atomic).
As a direct consequence of Corollary \ref{X-Min-discrete}, we obtain:

\begin{cor}
\label{Min-discrete}The following are equivalent for a commutative
(semi)ring $R:$

\begin{enumerate}
\item $R$ is an AMin-(semi)ring;

\item $Min(R)$ is discrete;

\item $P(R)=\bigcap\limits_{Q\in Min(R)}Q$ is an irredundant intersection.
\end{enumerate}
\end{cor}

While every FMin-(semi)ring is clearly an AMin-(semiring), the converse is
not true as the following example illustrates.

\begin{ex}
\label{AMin-Not-FMin}Let $\mathbb{F}$ be a field and $R:=\mathbb{F}%
[x_{1},x_{2},\cdots ]/(x_{i}x_{j}\mid i\neq j,$ $i,j\geq 1).$ Then
\begin{equation*}
Min(R)=\{P_{n}\mid n\geq 1\},\text{ where }P_{n}:=(\overline{x_{m}}\mid
m\neq n,m\geq 1).
\end{equation*}%
$P_{n}$ is absolutely minimal prime for every $n\geq 1$ since $\overline{x}%
_{n}\notin P_{n},$ while $\overline{x_{n}}\in \bigcap\limits_{m\neq n}P_{m}.$
So, $Min(R)=AMin(R)$ and $R$ is an AMin-ring which is not an FMin-ring.
\end{ex}

The following result follows by combining Theorem \ref{discrete} ad Remark %
\ref{Spec-discrete}:

\begin{thm}
\label{PAM}The following are equivalent for a commutative (semi)ring $R:$

\begin{enumerate}
\item $R$ is a PAMin-(semi)ring;

\item $R$ is a PBMax-(semi)ring;

\item $Spec(R)\ $is discrete;

\item $Spec(R)\ $is discrete and finite;

\item $Spec(R)\ $is finite and $T_{1}$ ($R_{0},$ $R_{1},$ $T_{2}$);

\item $R$ is an FMax-(semi)ring and $Spec(R)$ is $T_{1}$ ($R_{0},$ $R_{1},$ $%
T_{2}$);

\item $Max(R)\ $is discrete and $K.\dim (R)=0;$

\item $R$ is a BMax-(semi)ring and $K.\dim (R)=0$;

\item $R$ is an FMin-(semi)ring and $Spec(R)$ is $T_{1}$ ($R_{0},$ $R_{1},$ $%
T_{2}$);

\item $Min(R)$ is discrete and $K.\dim (R)=0;$

\item $R$ is an AMin-(semi)ring and $K.\dim (R)=0.$
\end{enumerate}
\end{thm}

\begin{cor}
\label{Spec(R)-discrete}The following are equivalent for a commutative ring $%
R:$

\begin{enumerate}
\item $Spec(R)\ $is discrete;

\item $Spec(R)\ $is discrete and finite;

\item $R$ is an FMax-ring and $\pi $-regular;

\item $R$ is an FMin-ring and $\pi $-regular.
\end{enumerate}
\end{cor}

The following example shows that the characterization of the discreteness of
$Spec(R)\ $in Corollary \ref{Spec(R)-discrete} may fail even for proper
semirings. It shows also that condition $K.\dim (R)=0$ cannot be dropped
(replaced by Artinian or $\pi $-regular) in the conditions that involve it
in Theorem \ref{PAM}.

\begin{ex}
\label{S3}(\cite{Abu2025}) Consider $(S,+,\cdot )$, where $S=\{0,a,1\}$ and
the addition and the multiplication are given by%
\begin{equation*}
\begin{tabular}{|l|l|l|l|}
\hline
$+$ & $0$ & $a$ & $1$ \\ \hline
$0$ & $0$ & $a$ & $1$ \\ \hline
$a$ & $a$ & $a$ & $1$ \\ \hline
$1$ & $1$ & $1$ & $1$ \\ \hline
\end{tabular}%
\text{ and }%
\begin{tabular}{|l|l|l|l|}
\hline
$\cdot $ & $0$ & $a$ & $1$ \\ \hline
$0$ & $0$ & $0$ & $0$ \\ \hline
$a$ & $0$ & $a$ & $a$ \\ \hline
$1$ & $0$ & $a$ & $1$ \\ \hline
\end{tabular}%
.
\end{equation*}%
Then $S$ is a commutative proper semidomain and%
\begin{equation*}
Ideal(S)=\{\{0\},\{0,a\},S\},\text{ }Spec(S)=\{\{0\},\{0,a\}\}\text{ and }%
Max(S)=\{\{0,a\}\}.
\end{equation*}

\begin{enumerate}
\item $S$ is subtractive and local.

\item $S$ is idempotent (whence von Neumann regular, $\pi $-regular).

\item $S$ is an FMin-semiring (whence AMin-semiring)\ and an FMax-semiring
(whence BMax-semiring).

\item $S$ is finite, whence Artinian and Noetherian.

\item $S$ is a semidomain (has no non-zero zero-divisors), In particular, $S$
is reduced.

\item $J(R)=\{0,a\}$ is \emph{not} a nil-ideal.

\item Both $Min(S)=\{\{0\}\}$ and $Max(S)=\{\{0,a\}\}$ have the discrete
topology.

\item The Zariski topology on $Spec(S)$ is $\tau _{Z}(S)=\{\emptyset ,$ $%
Spec(S),$ $\{0\}\}$ (\emph{not} discrete).

\item $K.\dim (S)=1$ ($S$ is \emph{not} Krull $0$-dimensional).

\item $Spec(S)$ is $T_{0}$ but not $T_{1}$ (notice that $\{0\}$ is not
closed).$\blacksquare $
\end{enumerate}
\end{ex}

\begin{defn}
We call a topological space $X$ \textbf{Noetherian }iff $\mathcal{O}(X)$
satisfies the Ascending Chain Condition (ACC); equivalently, $\mathcal{C}(X)$
satisfies the Descending Chain Condition (DCC).
\end{defn}

\begin{lem}
\label{X-Noeth}Let $\mathcal{L}=(L,\wedge ,\vee ,1,0)$ be an $X$-top lattice
for some $X\subseteq L\backslash \{1\}.$

\begin{enumerate}
\item $X$ is Noetherian if and only if the set $C^{X}(L)$ of $X$-radical
elements in $L$ satisfies the ACC.

\item If $X=SI^{C^{X}(L)}(C^{X}(L))$ and $X$ is Noetherian, then $Min(X)$ is
finite.

\item If $X$ is spectral and Noetherian, then $Min(X)$ is finite.
\end{enumerate}
\end{lem}

\begin{Beweis}
\begin{enumerate}
\item By \cite[Theorem 2.5 (1)]{AH2019}, we have a 1-1 correspondence
between the $X$\emph{-radical elements} of $L$ and the \emph{closed subsets}
of $X:$%
\begin{equation*}
V_{X}(-):C^{X}(L)\longrightarrow \mathcal{C}(X)\text{ and }I_{X}(-):\mathcal{%
C}(X)\longrightarrow C^{X}(L).
\end{equation*}

\item Recall that we denoted with $\mathcal{I}(X)$ (resp. $\mathcal{IC}(X)$)
the collection of irreducible closed subsets of $X$ (resp. the irreducible
components of $X$). If $X=SI^{C^{X}(L)}(C^{X}(L)),$ then it follows by (\cite%
[Theorem 2.5 (1)]{AH2019}) that the above 1-1 correspondence restricts to a
1-1 correspondences%
\begin{equation*}
X\longleftrightarrow \mathcal{I}(X)\text{ and }Min(X)\longleftrightarrow
\mathcal{IC}(X)
\end{equation*}%
Since $X$ is Noetherian, we conclude that $\left\vert Min(X)\right\vert
=\left\vert \mathcal{IC}(X)\right\vert <\infty $ (cf. \cite[Lemma 5.9.2 (Tag
0052)]{SPA} for the finiteness of $\mathcal{IC}(X)$ for any Noetherian
topological space $X$).

\item If $X$ is spectral, then $X=SI^{C^{X}(L)}(C^{X}(L))$ and the result
follow by (2).$\blacksquare $
\end{enumerate}
\end{Beweis}

The following result recovers \cite[Theorem 88]{Kap1974} and extends it to
commutative semirings.

\begin{cor}
\label{R-Noeth-Min}Let $R$ be a commutative (semi)ring.

\begin{enumerate}
\item $Spec(R)$ is Noetherian if and only if $R$ satisfies the ACC on its
radical ideals.

\item If $Spec(R)\ $is Noetherian, then $Min(R)$ is finite.

\item If $R$ is Noetherian, then $R$ is an FMin-(semi)ring.
\end{enumerate}
\end{cor}

\begin{thm}
\label{Artinian}The following are equivalent for a commutative \emph{%
Noetherian} ring $R:$

\begin{enumerate}
\item $R$ is Artinian;

\item $R$ is semilocal (i.e. $R/Jac(R)$ is semisimple)\ and $J(R)$ is
nilpotent;

\item $R$ is an FMax-ring and $J(R)$ is nilpotent;

\item $R$ is an FMax-ring and $J(R)$ is a nil ideal;

\item $K.\dim (R)=0;$

\item $Spec(R)$ is finite and $T_{1}$ ($R_{0},$ $R_{1},$ $T_{2}$);

\item $Spec(R)$ is discrete (and finite).

\item $R$ is a PAMin-ring;

\item $R$ is a PBMax-ring;

\item $R$ is $\pi $-regular.
\end{enumerate}
\end{thm}

\begin{Beweis}
Let $R$ be a Noetherian commutative ring.

$(1\Longleftrightarrow 2)$ An \emph{associative} ring $R$ is left\ Artinian
if and only if $R$ is left Noetherian, $R/J(R)$ is (left) semisimple and $%
J(R)\ $is nilpotent (e.g. \cite[31.4]{Wis1991}).

$(2\Longleftrightarrow 3)$ By \cite[Proposition 20.2]{Lam2001}, a \emph{%
commutative} ring semilocal if and only if $\left\vert Max(R)\right\vert
<\infty .$

$(3\Longrightarrow 4)$ is trivial.

$(4\Longrightarrow 5)\ $Assume that $Max(R)=\{\mathfrak{m}_{1},\cdots ,%
\mathfrak{m}_{n}\}$ and $J(R)=\bigcap\limits_{i=1}^{n}\mathfrak{m}_{i}$ is a
nil ideal. Then%
\begin{equation*}
\mathfrak{m}_{1}\cdot \cdots \cdot \mathfrak{m}_{n}\subseteq
\bigcap\limits_{i=1}^{n}\mathfrak{m}_{i}=J(R)\subseteq Nil(R)\overset{R\text{
is commutative}}{=}P(R).
\end{equation*}%
If $P$ is a prime ideal of $R,$ then it follows that $P=\mathfrak{m}_{i}$
for some $i.$ Consequently, $Spec(R)=Max(R),$ i.e. $R$ is Krull $0$%
-dimensional.

$(5\Longleftrightarrow 1)$ It's well known that a \emph{commutative} ring $R$
is Artinian if and only if $R$ is Noetherian and $K.\dim (R)=0$ (e.g. \cite[%
Theorem 8.5]{AM1969}).

$(5\Longleftrightarrow 6)$ Since $R$ is Krull $0$-dimensional, it turns out
the $\left\vert Spec(R)\right\vert =\left\vert Min(R)\right\vert <\infty $
since every commutative Noetherian (semi)ring has finitely many minimal
primes by Corollary \ref{R-Noeth-Min}. The result follows then by
Proposition \ref{dim-0}.

$(6\Longleftrightarrow 7\Longleftrightarrow 8\Longleftrightarrow
9\Longleftrightarrow 10)$ The equivalences follow from Theorem \ref{PAM}
taking into consideration that a Noetherian commutative ring is an FMin-ring
by Corollary \ref{R-Noeth-Min}.$\blacksquare $
\end{Beweis}

\begin{ex}
\label{Ex-Art-ring}Let $\mathbb{F}$ be a field and consider the Noetherian
ring%
\begin{equation*}
R:=\mathbb{F}[x]/(x^{2})\times \mathbb{F}[y]/(y^{3}).
\end{equation*}%
Notice that $R$ is an FMax-ring (semilocal) with
\begin{equation*}
Max(R)=\{\mathfrak{m}_{1},\mathfrak{m}_{2}\}=\{(\overline{x})\times \mathbb{F%
}[y]/(y^{3}),\text{ }\mathbb{F}[x]/(x^{2})\times (\overline{y})\}.
\end{equation*}%
Moreover, the \emph{Jacobson radical}%
\begin{equation*}
J(R)=J(\mathbb{F}[x]/(x^{2}))\times J(\mathbb{F}[y]/(y^{3}))=(\overline{x}%
)\times (\overline{y})
\end{equation*}%
is a \emph{nil-ideal} (one can see easily that $(a,b)^{3}=(\overline{0},%
\overline{0})$ for any $(a,b)=(\alpha \overline{x},\beta \overline{y}+\gamma
\overline{y}^{2})\in J(R)$). Notice that $R$ is Artinian as both $\mathbb{F}%
[x]/(x^{2})$ and $\mathbb{F}[y]/(y^{3})$ are \emph{finite dimensional}
algebras over $\mathbb{F}$.$\blacksquare $
\end{ex}

The following example shows that the if the assumption "$R$ \emph{is
Noetherian}" is removed from the context of Theorem \ref{Artinian}, then the
result may not be true by demonstrating the existence of a \emph{local}
commutative ring $R$ with $K.\dim (R)=0$ that is \emph{not} Artinian
(equivalently, $K.\dim (R)=0$ and $R$ is \emph{not} Artinian).

\begin{ex}
\label{AM-not-Art}(\cite[Page 91]{AM1969})\ Let $\mathbb{F}$ be a field, and
consider the commutative ring
\begin{equation*}
R:=\mathbb{F}[x_{1},x_{2},\cdots ,x_{n},\cdots ]/(x_{1},x_{2}^{2},\cdots
,x_{n}^{n},\cdots ).
\end{equation*}%
Then $Spec(R)=\{\mathfrak{m}\},$ where $\mathfrak{m}:=(\overline{x_{1}},%
\overline{x_{2}},\cdots ,\overline{x_{n}},\cdots ).$ So, $Spec(R)=Max(R)$ is
discrete and $(R,\mathfrak{m})$ is a local Krull $0$-dimensional ring.
Clearly, $\mathfrak{m}$ is not finitely generated, whence $R$ is \emph{not}
Noetherian (consequently,\ \emph{not} Artinian).$\blacksquare $
\end{ex}

By Theorem \ref{Artinian}, a $\pi $-regular Noetherian ring commutative $R$
is Krull $0$-dimensional and Artinian. The following example illustrates
that even a subtractive local idempotent (whence von Neumann regular, $\pi $%
-regular)\ commutative proper semiring can be neither Krull $0$-dimensional
nor Artinian.

\begin{ex}
\label{Noeth-not-Art-dim0}(\cite{Abu2025}) Consider $(S;\oplus ,\infty
;\otimes ,-\infty ),$ where%
\begin{equation*}
S=\mathbb{W}\cup \{-\infty ,\infty \},\text{ }a\oplus b:=\min \{a,b\}\text{
and }a\otimes b:=\max \{a,b\}.
\end{equation*}%
Then $S$ is a commutative proper semidomain with $0_{S}=\infty $ and $%
1_{S}=-\infty .$ Set
\begin{equation*}
G_{s}:=\{x\in S\mid x\geq s\}\text{ for }s\in S\text{.}
\end{equation*}%
Notice that $G_{-\infty }=S$ and $\mathbb{W}\cup \{\infty \}=G_{0}.$ It can
be shown that%
\begin{equation*}
\begin{tabular}{lllllll}
$Ideal(S)$ & $=$ & $\{G_{s}\text{ }\mid \text{ }s\in S\};$ &  & $\text{ }%
Max(S)$ & $=$ & $\{G_{0}\};\text{ }$ \\
$Spec(S)$ & $=$ & $\{G_{s}\text{ }\mid \text{ }s\geq 0\};$ &  & $Min(S)$ & $%
= $ & $\{G_{\infty }\}.$%
\end{tabular}%
\end{equation*}%
It follows that $S$ is subtractive and local. Moreover, $S$ is idempotent
(whence von Neumann regular, $\pi $-regular) and Noetherian but \emph{%
neither }Krull $0$-dimensional (in fact $K.\dim (S)=\infty $) \emph{nor}
Artinian: we have a strictly decreasing chain of (prime) ideals%
\begin{equation*}
G_{0}\varsupsetneqq G_{1}\varsupsetneqq \cdots \varsupsetneqq
G_{n}\varsupsetneqq G_{n+1}\varsupsetneqq \cdots .\blacksquare
\end{equation*}
\end{ex}

The key point in the proof of Theorem \ref{Artinian} is the well-known
characterization of commutative Artinian rings as the Krull $0$-dimensional
Noetherian ones. The following example illustrates that even a subtractive
local idempotent semidomain can be \emph{neither} Krull $0$-dimensional
\emph{nor} Noetherian.

\begin{ex}
\label{Art-not-Noeth}(\cite{Abu2025}) Consider $(S,\oplus ,0,\otimes ,\infty
),$ where%
\begin{equation*}
S=\mathbb{W}\cup \{\infty \},\text{ }a\oplus b:=\max \{a,b\}\text{ and }%
a\otimes b:=\min \{a,b\}.
\end{equation*}%
Then $S$ is a commutative proper semidomain with $0_{S}=0,$ $1_{S}=\infty .$
Set
\begin{equation*}
J_{s}:=\{x\in S\mid x\leq s\}\text{ for }s\in S.
\end{equation*}%
Notice that $S=J_{\infty }$ and $\mathbb{W}=\bigcup\limits_{n\in \mathbb{W}%
}J_{n}.$. It can be shown that%
\begin{equation*}
\begin{tabular}{lllllll}
$Ideal(S)$ & $=$ & $\{J_{n}\text{ }\mid \text{ }n\in \mathbb{W}\}\cup \{%
\mathbb{W},S\};$ &  & $Max(S)$ & $=$ & $\{\mathbb{W}\};$ \\
$\text{ }Spec(S)$ & $=$ & $\{J_{n}\text{ }\mid \text{ }n\in \mathbb{W}\}\cup
\{\mathbb{W}\};$ &  & $Min(S)$ & $=$ & $\{\{0\}\}.$%
\end{tabular}%
\end{equation*}%
It follows that $S$ is subtractive, local and idempotent (whence von Neumann
regular, $\pi $-regular). Moreover, $S$ is Artinian but \emph{neither Krull }%
$0$\emph{-dimensional }(in fact $K.\dim (S)=\infty $) \emph{nor Noetherian: }%
we have an infinite strictly ascending chain of (prime) ideals%
\begin{equation*}
J_{0}\subsetneqq J_{1}\subsetneqq \cdots \subsetneqq J_{n}\subsetneqq
J_{n+1}\subsetneqq \cdots .\blacksquare
\end{equation*}
\end{ex}

\section{Quarter Separation Axioms}

\qquad Let $\mathcal{L}=(L,\wedge ,\vee ,1,0)$ be an $X$-top lattice for
some $X\subseteq L\backslash \{1\}.$ By Proposition \ref{dim-0}, $X$ is $%
T_{0}$; however, $X$ is $T_{1}$ if and only if $K.\dim (X)=0.$ This suggests
studying the separation axioms between $T_{0}$ and $T_{1}.$ While there are
many such separation axioms (cf. \cite{Avi2005}, \cite{Avi2006}), we limit
our attention in this paper to the so-called \textbf{quarter separation
axioms}: $T_{\frac{1}{4}},$ $T_{\frac{1}{2}}$ and $T_{\frac{3}{4}}$ as
defined in \cite{PRV2009}.

\bigskip

For the convenience of the reader, we start by recalling and fixing the
needed definitions.

\begin{defn}
(e.g. \cite{PRV2009}, \cite{Dun1977}) We say that a topological space $X$ is

\begin{enumerate}
\item $T_{\frac{1}{4}}$ iff any $x\in X$ is closed \emph{or} kerneled (i.e.
iff $X=Cl(X)\cup K(X)$).

\item $T_{\frac{1}{2}}$ iff any $x\in X$ is closed \emph{or} isolated (i.e.
iff $X=Cl(X)\cup Iso(X)$).

\item $T_{\frac{3}{4}}$ iff any $x\in X$ is closed \emph{or} regular open
(i.e. iff $X=Cl(X)\cup RO(X)$).
\end{enumerate}
\end{defn}

\begin{rem}
(e.g.$\ $\cite{PRV2009}) Notice that for any topological space $X,$ we have $%
RO(X)\subseteq Iso(X)\subseteq K(X),$ whence%
\begin{equation*}
T_{1}\ \implies T_{\frac{3}{4}}\ \implies T_{\frac{1}{2}}\ \implies T_{\frac{%
1}{4}}\ \implies T_{0}.
\end{equation*}
\end{rem}

\bigskip

\begin{defn}
(cf. \cite{Avi2005}) Let $X$ be a topological space. For $A,B\subseteq X,$
we set $A\vdash B$ iff there exists $U\underset{\text{open}}{\subseteq }X$
such that $A\subseteq U$ and $B\cap U=\emptyset .$ We say that $X$ is $T_{F}$
iff for every $x\in X$ and $F\underset{\text{finite}}{\subseteq }X\backslash
\{x\},$ either $\{x\}\vdash F$ or $F\vdash \{x\}.$
\end{defn}

\bigskip

We begin by an easy, but very useful, diagrammatic characterization of $X$%
-top lattices for which $X$ is $T_{\frac{1}{4}}.$

\begin{thm}
\label{t1/4t}Let $L$ be an $X$-top lattice for some $X\subseteq L\backslash
\{1\}.$ The following are equivalent:

\begin{enumerate}
\item $X$ is a $T_{\frac{1}{4}};$

\item $K.\dim (X)\leq 1;$

\item $X$ is $T_{F}.$
\end{enumerate}
\end{thm}

\begin{Beweis}
$(1\Longleftrightarrow 2)$ We have%
\begin{equation*}
X\text{ is }T_{\frac{1}{4}}\Longleftrightarrow X\overset{\text{Definition}}{=%
}Cl(X)\cup K(X)\overset{\text{Lemma \ref{kerneled}}}{=}Max(X)\cup
Min(X)\Longleftrightarrow K.\dim (X)\leq 1.
\end{equation*}

$(2\Longrightarrow 3)$ Assume that $K.\dim (X)\leq 1.$ Let $x\in X$ and $F%
\underset{\text{finite}}{\subseteq }X\backslash \{x\}.$ If $F=\emptyset ,$
then $\{x\}\vdash F$ as $x\in X$ and $F\cap X=\emptyset .$ Assume that $%
F\neq \emptyset $ and $F:=\{x_{1},\cdots ,x_{n}\}\subseteq X\backslash
\{x\}. $

\textbf{Case I:}\ $x\in Min(X).$ In this case, $x\in
D(\bigwedge\limits_{i=1}^{n}x_{i})$ (since $X$ is strongly $X$-irreducible
in $C^{X}(L)$ by Theorem \ref{xct} and $X\subseteq C^{X}(L)$). It's obvious
that $D(\bigwedge\limits_{i=1}^{n}x_{i})\cap F=\emptyset .$ So, $\{x\}\vdash
F.$

\textbf{Case II:\ }$x\in Max(X).$ In this case, $V(x)=\{x\}\ $whence $%
F\subseteq X\backslash \{x\}=D(x)$ and $\{x\}\cap D(x)=\emptyset ,$ i.e. $%
F\vdash \{x\}.$

$(3\Longrightarrow 1)$ Assume that $X$ is $T_{F}.$ Suppose that there exists
$x\in X\backslash (Max(X)\cup Min(X)).$ Then there exists $%
F:=\{x_{1},x_{2}\}\subseteq X\backslash \{x\}$ such that $x_{1}\lneqq
x\lneqq x_{2}.$

\textbf{Case I:}$\ \{x\}\vdash F.$ Then there exists $a\in L$ such $x\in
D(a) $ and $F\cap D(a)=\emptyset .$ It follows that $a\leq x_{1}\lneqq x,$ a
contradiction to $x\in D(a).$

\textbf{Case II:}$\ F\vdash \{x\}.$ Then there exists $b\in L$ such that $%
F\subseteq D(b)$ and $x\notin D(b).$ It follows that $b\leq x\lneqq x_{2},$
a contradiction to $x_{2}\in D(b).\blacksquare $
\end{Beweis}

\bigskip

The following consequence of Theorem \ref{t1/4t}, recovers \cite[Corollary
5.1]{PRV2009} and extends \cite[Proposition 3.1.3]{Avi2005} to commutative
semirings.

\bigskip

\begin{cor}
\label{semi-TF}Let $R$ be a commutative (semi)ring. Then%
\begin{equation*}
Spec(R)\text{ is }T_{\frac{1}{4}}\Longleftrightarrow K.\dim (R)\leq
1\Longleftrightarrow Spec(R)\text{ is }T_{F}.
\end{equation*}
\end{cor}

The following example shows that a $T_{\frac{1}{4}}$ spectral space is not
necessarily $T_{1}.$

\begin{ex}
\label{Z-anti}Consider the ring $\mathbb{Z}$ of integers, $L=Ideal(\mathbb{Z}%
),$ $X=Spec(\mathbb{Z})$ and $Y:=Max(\mathbb{Z}):=Spec(\mathbb{Z})\backslash
\{0\}.$Since $K.\dim (Spec(\mathbb{Z}))=1,$ we conclude by that $Spec(%
\mathbb{Z})$ is $T_{\frac{1}{4}}$ but not $T_{1}$ (by Theorem \ref{t1/4t}
and Proposition \ref{dim-0} (4)). Moreover, $\mathbb{Z}$ is a $Max(\mathbb{Z}%
)$\emph{-top ring} by Corollary \ref{Y-top}. Since $K.\dim (Max(\mathbb{Z}%
))=0,$ we conclude, by Proposition \ref{dim-0}, that $Max(\mathbb{Z})$ is $%
T_{1}.$ Notice that $Max(\mathbb{Z})$ is not $T_{2},$ whence not spectral
since every $T_{1}$ spectral space is $T_{2}$ by Lemma \ref{T1-qH-T2}. In
fact, $Max(\mathbb{Z})$ is far away from being $T_{2}$ as it is anti-$T_{2}.$

\textbf{Claim:}\ $Y:=Max(\mathbb{Z})$ is hyperconnected (irreducible).

Let $D(n\mathbb{Z})$ and $D(m\mathbb{Z})$ be two \emph{disjoint} open sets
(for some positive integers $m\neq n$). Then%
\begin{equation*}
\begin{array}{ccccc}
Y\backslash V(nm\mathbb{Z}) & = & Y\backslash (V(n\mathbb{Z}\cap m\mathbb{Z}%
)) & = & Y\backslash (V(n\mathbb{Z})\cup V(m\mathbb{Z})) \\
& = & Y\backslash V(n\mathbb{Z})\cap Y\backslash V(m\mathbb{Z}) & = &
\emptyset .%
\end{array}%
\end{equation*}%
Hence, $V(nm\mathbb{Z})=Y$, i.e. every prime number divides $nm$ which is
absurd. Therefore, $Y$ is anti-$T_{2}$.$\blacksquare $
\end{ex}

\begin{ex}
\label{W-anti}Let $L$ be any $X$-top lattice with $K.\dim (X)>1.$ Then $X$
is $T_{0}$ but not $T_{\frac{1}{4}}$ (by Proposition \ref{dim-0} and Theorem %
\ref{t1/4t}). Any subspace of $X$ is indeed $T_{0}$ and could be $T_{\frac{1%
}{4}}$ or even $T_{1}.$

Consider prime spectrum of the commutative semiring $\mathbb{W}$ of whole
numbers shown in Figure \ref{spec-w}.

\begin{figure}[tbp]
\centering
\begin{equation*}
\begin{tabular}{c}
$\xymatrix{& & & {\mathbb{W}}\backslash{{\{1\}}} \ar@{-}[ddll] \ar@{-}[ddl]
\ar@{-}[dd] \ar@{-}[ddr] \ar@{-}[ddr] \ar@{-}[ddlll] & \\ & & & & \\
2{\mathbb{W}} \ar@{-}[ddrrr] & 3{\mathbb{W}} \ar@{-}[ddrr] & \cdots
\ar@{-}[ddr] & p {\mathbb{W}} \ar@{-}[dd] & \cdots \ar@{-}[ddl] \\ & & & &
\\ & & & 0 & }$ \\
\\
\end{tabular}%
\end{equation*}%
\caption{\textbf{The prime spectrum of} $\mathbb{W}$}
\label{spec-w}
\end{figure}
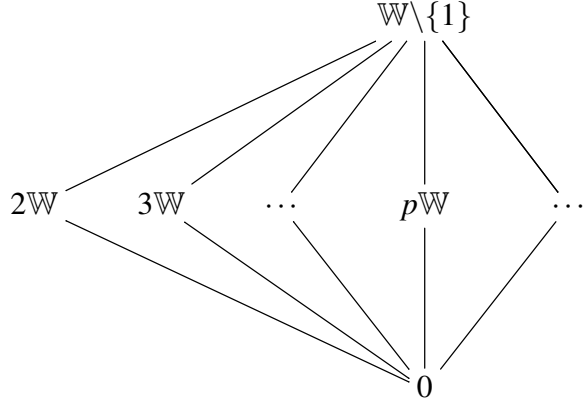

Since $K.\dim (\mathbb{W})=2,$ we conclude that $Spec(\mathbb{W})$ is $T_{0}$
but not $T_{\frac{1}{4}}.$

Set $Y:=Spec(\mathbb{W})\backslash \{\mathbb{W}\backslash \{1\}\}.$ It
follows by Corollary \ref{Y-top} that $\mathbb{W}$ is a $Y$-\emph{top
semiring}. Moreover, $K.\dim (Y)=1,$ whence $Y$ is $T_{\frac{1}{4}}$ but not
$T_{1}.$

Set $Z:=Spec(\mathbb{W})\backslash \{\mathbb{W}\backslash \{1\},0\}.$ It
follows by Corollary \ref{Y-top} that $\mathbb{W}$ is a $Z$-\emph{top
semiring}. Moreover, $K.\dim (Z)=0,$ whence $Z$ is $T_{1}.\blacksquare $
\end{ex}

\begin{defn}
Let $L$ be an $X$-top lattice for some $X\subseteq L\backslash \{1\}.$ We
say that $X$ is $ES$ iff $Min(X)\backslash Max(X)\subseteq CSI^{X}(X),$ i.e.
for every $m\in Min(X)\backslash Max(X),$ we have%
\begin{equation*}
\bigwedge\limits_{x\in X\backslash \{m\}}x\nleqq m.
\end{equation*}%
We call a (semi)ring $R$ an $ES$-(semi)ring iff $Spec(R)$ is $ES$ (\cite%
{Avi2006}).
\end{defn}

Next, we investigate the $X$-top lattices for which $X$ is $T_{\frac{1}{2}}$.

\begin{thm}
\label{T1/2}Let $L$ be an $X$-top lattice for some $X\subseteq L\backslash
\{1\}.$ The following are equivalent:

\begin{enumerate}
\item $X$ is $T_{\frac{1}{2}};$

\item $X=Max(X)\cup (Min(X)\cap CSI^{X}(X)).$

\item $K.\dim (X)\leq 1$ and $X$ is $ES.$

\item $X$ is $T_{\frac{1}{4}}$ and $X$ is $ES.$
\end{enumerate}
\end{thm}

\begin{Beweis}
$(1\Longleftrightarrow 2)$ It follows from the definition, Lemma \ref%
{kerneled} and Proposition \ref{iso} that%
\begin{equation*}
X\text{ is }T_{\frac{1}{2}}\Longleftrightarrow X=Cl(X)\cup Iso(X)=Max(X)\cup
(Min(X)\cap CSI^{X}(X)).
\end{equation*}

$(2\Longleftrightarrow 3)$ This is obvious.

$(3\Longleftrightarrow 4)$ This follows from Theorem \ref{t1/4t}.$%
\blacksquare $
\end{Beweis}

\begin{ex}
\label{rai}Let $L$ be an $X$-top lattice for some $X\underset{\emph{finite}}{%
\subseteq }L\backslash \{1\}$ and $X$ is . Since $X$ is finite, we have $%
X=CSI^{X}(X)$ by Remark \ref{finite-BMax} (2), whence $X$ is, in particular,
$ES.$
\end{ex}

\begin{cor}
\label{T1/4-T1/2}Let $L$ be an $X$-top lattice for some $X\subseteq
L\backslash \{1\}.$ If $X$ is $ES$ (e.g. $X$ is \emph{finite}), then
\begin{equation*}
X\text{ is }T_{\frac{1}{2}}\Longleftrightarrow K.\dim (X)\leq
1\Longleftrightarrow X\text{ is }T_{\frac{1}{4}}\Longleftrightarrow X\text{
is }T_{F}.
\end{equation*}
\end{cor}

The following result recovers \cite[Proposition 5]{Avi2006} and extends it
to commutative semirings:

\begin{cor}
\label{ES-semiring}Let $R$ be a commutative (semi)ring. Then $Spec(R)$ is $%
T_{\frac{1}{2}}$ if and only if $R$ is an $ES$-(semi)ring and $K.\dim
(R)\leq 1.$
\end{cor}

The conditions: $R$ is an $ES$ semi(ring) and $K.\dim (R)\leq 1$ in
Corollary \ref{ES-semiring} are \emph{independent} as the following two
examples show, whence none of them can be dropped.

\begin{ex}
\label{valuation}Let $D$ be a \emph{valuation ring} (i.e. an integral domain
with field of fractions $K$ such that for every $x\in K\backslash \{0\}$
either $x\in D$ or $x^{-1}\in D$ (equivalently, $(Ideal(D),\subseteq )$ is
totally ordered \cite[Chapter VI. (1.2.), Theorem 1]{Bou1989}). It is clear
that $Min(D)=\{0\}$ and $\bigcap\limits_{P\in Spec(D)\backslash
\{0\}}P\nsubseteqq 0,$ i.e. $D$ is $ES.$ So, $Spec(D)$ is $T_{\frac{1}{2}}$
if and only if $K.\dim (D)\leq 1.$ So, if $D$ is any valuation ring with $%
K.\dim (D)\geq 2,$ then $D$ is $ES$ but $Spec(D)$ is not $T_{\frac{1}{2}}.$
\end{ex}

\begin{ex}
\label{DVR-valuation}A DVR (\textbf{discrete valuation ring}) is a local PID
that is \emph{not} a field, equivalently a Noetherian local valuation ring
that is distinct from its field of fractions \cite[Chapter VI, 3.6,
Proposition 9]{Bou1989}. In particular, every DVR is an $ES$-ring. If $(D,%
\mathfrak{m})$ is a DVR, then $Spec(D)=\{0,\mathfrak{m}\}$ by \cite[page 94]%
{AM1969}, whence $Spec(D)$ is $T_{\frac{1}{2}}$ by Corollary \ref%
{ES-semiring}. One can check this easily that the Zariski topology on $%
X=Spec(D)$ is $\tau =\{\emptyset ,X,\{0\}\},$ i.e. $\mathfrak{m}$ is closed
and $0$ is isolated. A concrete example of a DVR is $D=\mathbb{F}[[x]]$: the
ring of formal power series with coefficients in a field $\mathbb{F}$. For
more examples of DVRs, the reader may consult \cite[Chapter VI. (1.2.),
Theorem 1]{Bou1989} (and \cite{DF2004}).
\end{ex}

\begin{ex}
$Spec(\mathbb{Z})$ is $T_{\frac{1}{4}}$ (as $K.\dim (\mathbb{Z})=1$)\ but
not $T_{\frac{1}{2}}.$ In fact, $\{0\}\neq D((m))$ for any $m\in \mathbb{Z},$
whence $\{0\}\notin Iso(Spec(\mathbb{Z})).$ Notice that $Min(\mathbb{Z}%
)=\{0\}$ but $\bigcap\limits_{p\ \text{prime}}p\mathbb{Z}=0,$ whence $Spec(%
\mathbb{Z})$ is \emph{not} $ES$ (cf. Theorem \ref{T1/2}). One can show
similarly (and for the same reason) that $Y=Spec(\mathbb{W})\backslash (%
\mathbb{W}\backslash \{1\})$ is $T_{\frac{1}{4}}$ but not $T_{\frac{1}{2}}$ (%
$Y$ is \emph{not} $ES$).
\end{ex}

Before we proceed, we recall some definitions from \emph{Graph Theory} that
will be used throughout the rest of the section.

\begin{punto}
Let $(P,\leq )$ be a partially ordered set. We call $x_{0}<x_{1}<\cdots
<x_{n}$ in $P$ a \textbf{chain of length }$n$ and denote it with $\mathcal{C}%
_{n+1}.$ We consider also spacial types of \textbf{trees with finite base} $%
\mathcal{T}_{n}$ consisting of one maximal element lying over $n$ minimal
elements\ and \textbf{dual trees with finite covers} $\mathcal{V}_{n}$
consisting of one minimal element lying under $n$-maximal elements (a tree $%
\mathcal{T}_{n}$ in the dual poset $(P,\geq )$).

\begin{figure}[tbp]
\begin{equation*}
\begin{array}{ccc}
\xymatrix{& & \bullet \ar@{-}[dll] \ar@{-}[dl] \ar@{-}[d] \ar@{-}[dr]
\ar@{-}[drr] & & \\ \bullet & \bullet & \bullet & \bullet & \bullet} & %
\xymatrix{\bullet \ar@{-}[d] \\ \bullet } & \xymatrix{\bullet & \bullet &
\bullet \\ & \bullet \ar@{-}[ul] \ar@{-}[u] \ar@{-}[ur] } \\
\mathcal{T}_{5} & \mathcal{T}_{1}=\mathcal{C}_{2}=\mathcal{V}_{1} & \mathcal{%
V}_{3}%
\end{array}%
\end{equation*}%
\caption{\textbf{Examples of trees and dual trees}}
\label{trees-dual}
\end{figure}
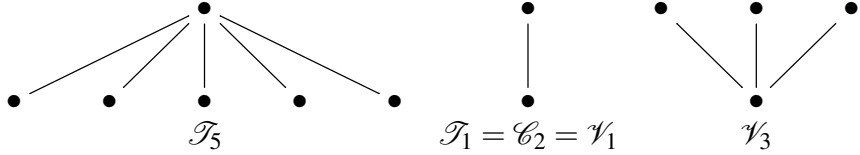
Examples are given in Figure \ref{trees-dual}. A collection of \emph{disjoint%
} trees is called a \textbf{forest}.
\end{punto}

The following result provides a big family of $X$-top lattices for which $X$
is $T_{\frac{1}{2}}$ but \emph{not} $T_{1}.$

\begin{ex}
\label{TV}Let $L$ be an $X$-top lattice for some $X\subseteq L\backslash
\{1\}.$ Assume that $X$ is a \emph{finite} combination of \emph{disjoint} $%
\mathcal{T}_{n}$'s and $\mathcal{V}_{k}$'s. Notice that $X$ is a finite and $%
K.\dim (X)=1,$ So, $X$ is $T_{\frac{1}{4}}$ by Theorem \ref{t1/4t}, whence $%
T_{\frac{1}{2}}$ by Corollary \ref{T1/4-T1/2}. Notice that $X$ is not $T_{1}$
by Proposition \ref{dim-0}.$\blacksquare $
\end{ex}

The following result extends partially \cite[Theorem 5.2]{PRV2009} from the
special case of the lattice of ideals of a commutative semiring to arbitrary
$X$-top lattices.

\begin{prop}
\label{RO-ISO-D}Let $L$ be an $X$-top lattice for some $X\subseteq
L\backslash \{1\}.$ Then%
\begin{equation*}
RO(X)=\{x\in X\mid \overline{D(x)}=X\backslash \{x\}\}=Iso(X)\cap \mathcal{E}%
(X).
\end{equation*}
\end{prop}

\begin{Beweis}
\textbf{Step I:}$\ RO(X)\subseteq \{x\in X\mid \overline{D(x)}=X\backslash
\{x\}\}.$

Let $x\in RO(X),$ i.e. $\{x\}=int(V(x)).$ Since $\{x\}\underset{\text{open}}{%
\subseteq }X$ and $D(x)\subseteq X\backslash \{x\},$ it follows that $%
\overline{D(x)}\subseteq X\backslash \{x\}.$ On the other hand, let $y\in
X\backslash \{x\}.$ If $U_{y}\in \mathcal{O}(y)$ and $U_{y}\subseteq V(x),$
then $y\in int(V(x))=\{x\},$ a contradiction. It follows that $U_{y}\cap
D(x)\neq \emptyset $ for every $U_{y}\in \mathcal{O}(y),$ whence $y\in
\overline{D(x)}.$ So, $X\backslash \{x\}\subseteq \overline{D(x)}.$
Consequently, $X\backslash \{x\}=\overline{D(x)}.$

\textbf{Step II:}$\ \{x\in X\mid \overline{D(x)}=X\backslash
\{x\}\}\subseteq Iso(X)\cap \mathcal{E}(X).$

Let $x\in X$ be such that $\overline{D(x)}=X\backslash \{x\}.$ Notice that $%
\{x\}=X\backslash \overline{D(x)}\underset{\text{open}}{\subseteq }X,$ i.e. $%
x\in Iso(X).$ Since $D(x)\subseteq X\backslash \{x\},$ we have $%
E(x):=\bigwedge (X\backslash \{x\})\leq \bigwedge D(x).$ On the other hand,
if $a_{x}:=\bigwedge D(x)\nleqq E(x),$ then there exists $q\in X\backslash
\{x\}\ $such that $a_{x}\nleqq q$, i.e. $q\in D(a_{x})$ and $D(a_{x})\in
\mathcal{O}(q).$ By assumption, $q\in X\backslash \{x\}=\overline{D(x)},$
whence $D(a_{x})\cap D(x)\neq \emptyset ,$ a contradiction.

\textbf{Step III:} $Iso(X)\cap \mathcal{E}(X)\subseteq RO(X).$

Let $x\in Iso(X)\cap \mathcal{E}(X).$ By assumption, $\{x\}\underset{\text{%
open}}{\subseteq }X,$ whence $x\in int(\overline{\{x\}}).$ On the other
hand, let $q\in int(\overline{\{x\}})=int(V(x)).$ Then $q\in D(a)\subseteq
V(x)$ for some $a\in L.$ If $q\neq x,$ then $q\in X\backslash \{x\}$ and so $%
q\in V(E(x))=V(\bigwedge D(x))=\overline{D(x)}.$ Since $D(a)\in \mathcal{O}%
(q),$ there exists $h\in D(a)\cap D(x)\subseteq V(x)\cap D(x)=\emptyset ,$ a
contradiction. So, $\{x\}=int(\overline{\{x\}}),$ i.e. $x\in
RO(X).\blacksquare $
\end{Beweis}

In what follows, we provide sufficient/necessary conditions for an $X$-top
lattice to be $T_{\frac{3}{4}}.$

The following characterization follows directly from the definition, Lemma %
\ref{kerneled} Proposition \ref{RO-ISO-D} and Proposition \ref{iso}:

\begin{thm}
\label{3/4T}Let $L$ be an $X$-top lattice for some $X\subseteq L\backslash
\{1\}.$ The following are equivalent:

\begin{enumerate}
\item $X$ is $T_{\frac{3}{4}};$

\item $X=Max(X)\cup RO(X);$

\item $X=Max(X)\cup (Iso(X)\cap \mathcal{E}(X));$

\item $X=Max(X)\cup (Min(X)\cap CSI^{X}(X)\cap \mathcal{E}(X)).$
\end{enumerate}
\end{thm}

Now, we provide a sufficient/necessary condition for an $X$-top lattice with
$X$ a \emph{finite forest }so that $X$ is (is not) $T_{\frac{3}{4}}.$

\begin{thm}
\label{t3/4t}Let $L$ be an $X$-top lattice for some $X\underset{\text{finite}%
}{\subseteq }L\backslash \{1\}.$

\begin{enumerate}
\item If $X$ is a forest consisting of a combination of \emph{disjoint} $%
\mathcal{T}_{n_{i}}$'s with $n_{i}\geq 2$ for each $i,$ then $X$ is $T_{%
\frac{3}{4}}$ but \emph{not} $T_{1}.$

\item If $X$ is $T_{\frac{3}{4}},$ then $K.\dim (X)\leq 1$ and $X$ does
\emph{not} contain a \emph{disjoint} $\mathcal{V}_{n}$.
\end{enumerate}
\end{thm}

\begin{Beweis}
\begin{enumerate}
\item Let $X=\bigcup\limits_{i=1}^{k}\mathcal{T}_{n_{i}}$ and let $\mathfrak{%
m}_{i}$ be the unique maximal element of $\mathcal{T}_{n_{i}}$ for $%
i=1,\cdots ,k.$ Then $Max(X)=\{\mathfrak{m}_{1},\cdots ,\mathfrak{m}_{n}\}$
and $Min(X)=\bigcup\limits_{i=1}^{k}(\mathcal{T}_{n_{i}}\backslash \{%
\mathfrak{m}_{i}\}).$ Let $x\in Min(X)$. Then there exists some $j\in
\{1,\cdots ,k\}$ such that $x\in \mathcal{T}_{n_{j}}\backslash \{\mathfrak{m}%
_{j}\}.$

\textbf{Claim 1:} $V(x)$ is \emph{not} open.

Suppose that $V(x)=\{x,\mathfrak{m}_{j}\}=D(z)$ for some $z\in L.$ Then $%
z\nleq \mathfrak{m}_{j},$ whence $\mathcal{T}_{n_{j}}\subseteq D(z).$ Since $%
n_{j}\geq 2,$ we have $D(z)=V(x)\varsubsetneqq \mathcal{T}_{n_{j}}\subseteq
D(z),$ a contradiction.

\textbf{Claim 2: }$\{x\}$ is open. This follows from%
\begin{equation*}
\{x\}=X\backslash (\bigcup\limits_{a\in X\backslash
\{x\}}V(a))=\bigcap\limits_{a\in X\backslash \{x\}}D(a)
\end{equation*}%
and the assumption that $X$ is finite. It follows that%
\begin{equation*}
int(\overline{\{x\}})=int(V(x))=int(\{x,\mathfrak{m}_{j}\})=\{x\},
\end{equation*}

i.e. $x$ is regular open. Consequently, $Min(X)=RO(X).$

So, $X=Max(X)\cup RO(X)$, i.e. $X$ is $T_{\frac{3}{4}}$ by Theorem \ref{3/4T}%
. Since $K.\dim (X)=1,$ it follows by Proposition \ref{dim-0} that $X$ is
\emph{not} $T_{1}.$

\item Let $X$ be $T_{\frac{3}{4}}.$ Then $X$ is $T_{\frac{1}{4}},$ whence $%
K.\dim (X)\leq 1$ by Theorem \ref{t1/4t}.

\begin{figure}[tbp]
\begin{equation*}
\begin{tabular}{ccc}
& $\xymatrix{ & & {\bullet} m_1 \ar@{-}[ddrrr] & {\bullet} m_2 \ar@{-}[ddrr]
& & & \cdots & {\bullet} m_n \ar@{-}[ddll] & \\ & & & & & & & & \\ & & & & &
\underset{x}{\bullet}}$ &
\end{tabular}%
\end{equation*}%
\caption{\textbf{A dual tree }$\mathcal{V}_{n}$}
\label{dual-tree}
\end{figure}
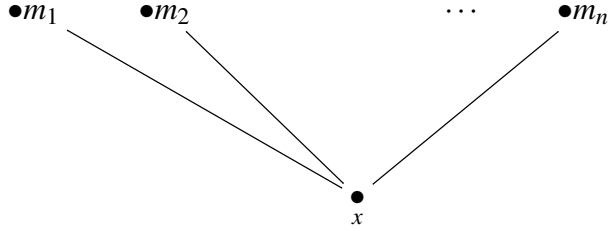

Suppose that $X$ contains a separate $\mathcal{V}_{n}$ with minimal element $%
x$ (notice that $\mathcal{V}_{1}=\mathcal{C}_{2}$) as shown in Figure \ref%
{dual-tree}.

Set $a:=\bigwedge\limits_{x\in X\backslash \mathcal{V}_{n}}.$ Then $a\nleq m$
for all $m\in Max(\mathcal{V}_{n})$ as $X=SI^{X}(X)$ by Remark \ref%
{finite-BMax} (1). It follows that $a\nleq x$ (the minimal element of $%
\mathcal{V}_{n}\subseteq X$). Therefore, $V(x)=\mathcal{V}_{n}=D(a)$ and $%
x\notin RO(X)$ since%
\begin{equation*}
int(\overline{\{x\}})=int(V(x))=V(x)=\mathcal{V}_{n}\neq \{x\}.
\end{equation*}%
Since $x$ is \emph{neither} maximal \emph{nor} regular open, it follows that
$X$ is \emph{not} $T_{\frac{3}{4}}$ (by Theorem \ref{3/4T}).$\blacksquare $
\end{enumerate}
\end{Beweis}

\begin{lem}
\label{localization-AM}\emph{(cf. \cite[Exercise 1.22 $\&\ $Proposition 3.11]%
{AM1969})}\ Let $R$ be a ring.

\begin{enumerate}
\item Let $R=\prod\limits_{i=1}^{n}R_{i}$ (a finite direct product of
rings). Then $Spec(R)\approx \bigsqcup\limits_{i=1}^{n}Spec(R_{i})$ (i.e. $%
Spec(R)$ is homeomorphic to the disjoint union of the prime spectra of the
rings $R_{1},\cdots ,R_{n}$).

\item If $S\subseteq R$ is a multiplicatively closed set, then the is a 1-1
correspondence%
\begin{equation*}
\{P\in Spec(R)\text{ }|\text{ }P\cap S=\emptyset \}\longleftrightarrow
Spec(S^{-1}R),\text{ }P\longmapsto S^{-1}P.
\end{equation*}
\end{enumerate}
\end{lem}

\begin{punto}
\label{Rn}Let $\mathbb{F}$ be a field and $n\geq 2.$ Consider%
\begin{equation*}
A_{n}:=\mathbb{F}[x_{1},\cdots ,x_{n}]/(x_{i}x_{j}\mid i\neq j),\text{ }%
\mathfrak{m}_{n}:=(\overline{x_{1}},\cdots ,\overline{x_{n}})\text{ and }%
R_{[n]}:=(A_{n})_{\mathfrak{m}}.
\end{equation*}%
For each $i,$ set $P_{i}:=(\overline{x_{1}},\cdots ,\overline{x_{i-1}},%
\overline{x_{i+1}},\cdots ,\overline{x_{n}}).$ Then $A_{n}/P_{i}\simeq
\mathbb{F}[x_{i}]$ is an integral domain and $A_{n}/\mathfrak{m}_{n}\simeq
\mathbb{F}$, hence $P_{i}$ is a (minimal) prime of $A_{n}$ ad $\mathfrak{m}%
_{n}$ is a maximal ideal of $A_{n}.$ The localization of $A_{n}$ at its
maximal ideal $\mathfrak{m}_{n}$ is a local ring $(R_{[n]},\mathfrak{m}%
_{n}R_{[n]}).$ Since $R_{[n]}/P_{i}R_{[n]}\simeq (\mathbb{F}%
[x_{i}])_{(x_{i})}$ and $Spec((\mathbb{F}[x_{i}])_{(x_{i})})=\{\{0\},(x_{i})%
\mathbb{F}[x_{i}])_{(x_{i})}\}$ for $i=1,\cdots ,n,$ it follows by Lemma \ref%
{localization-AM} (2) that%
\begin{equation*}
Spec(R_{[n]})=\{P_{[1]}R_{[n]},\cdots ,P_{[n]}R_{[n]},\mathfrak{m}%
_{n}R_{[n]}\}\text{ and }Max(R_{[n]})=\{\mathfrak{m}_{n}R_{[n]}\}.
\end{equation*}%
In particular, $Spec(R_{[n]})$ is a tree $\mathcal{T}_{n}.$
\end{punto}

\begin{punto}
\label{Pn}Let $\{p_{1},\cdots ,p_{n}\}$ be a set of distinct primes of $%
\mathbb{Z}$. Set%
\begin{equation*}
S:=\{m\in \mathbb{Z}\mid \gcd (m,p_{1}p_{2}\cdots p_{n})=1\}\text{ and }%
R:=S^{-1}\mathbb{Z}\text{.}
\end{equation*}%
It follows from Lemma \ref{localization-AM} (2) that $Spec(R)$ corresponds
to the primes of $\mathbb{Z}$ disjoint from $S$, i.e. $Spec(R)=\{\{0_{R}%
\},(p_{1})R,\cdots ,(p_{n})R\}$ with $(p_{i})R$ and $(p_{j})R$ incomparable
for $i\neq j.$ In particular, $Spec(R)\ $is a $\mathcal{V}_{n}.$
\end{punto}

\begin{exs}
\label{R2-R3}Let $\mathbb{F}$ be a field. Consider%
\begin{equation*}
R:=R_{[2]}\times R_{[3]}=(\mathbb{F}[x_{1},x_{2}]/(x_{1}x_{2}))_{(\overline{%
x_{1}},\overline{x_{2}})}\times (\mathbb{F}%
[y_{1},y_{2},y_{3}]/(y_{1}y_{2},y_{1}y_{3},y_{2}y_{3}))_{(\overline{y_{1}},%
\overline{y_{2}},\overline{y_{3}})}.
\end{equation*}%
Recall from \ref{Rn} that%
\begin{eqnarray*}
Spec(R_{[2]}) &=&\{(\overline{x_{2}})R_{[2]},(\overline{x_{1}})R_{[2]},(%
\overline{x_{1}},\overline{x_{2}})R_{[2]}\}; \\
Spec(R_{[3]}) &=&\{(\overline{y_{2}},\overline{y_{3}})R_{[3]},(\overline{%
y_{1}},\overline{y_{3}})R_{[3]},(\overline{y_{1}},\overline{y_{2}})R_{[3]},(%
\overline{y_{1}},\overline{y_{2}},\overline{y_{3}})R_{[3]}\}.
\end{eqnarray*}%
It follows by Lemma \ref{localization-AM} (1) that
\begin{equation*}
Spec(R)\approx Spec(R_{[2]})\bigsqcup Spec(R_{[3]}),
\end{equation*}%
a disjoint union $\mathcal{T}_{2}\bigsqcup \mathcal{T}_{3}.$ By Theorem \ref%
{t3/4t} (1), $Spec(R)$ is $T_{\frac{3}{4}}$ but not $T_{1}.\blacksquare $
\end{exs}

\begin{exs}
\label{S1-S2}Consider
\begin{eqnarray*}
S_{1} &:&=\{m\in \mathbb{Z}\mid \gcd (m,6)=1\},\text{ }R_{1}:=S_{1}^{-1}%
\mathbb{Z}; \\
S_{2} &:&=\{m\in \mathbb{Z}\mid \gcd (m,70)=1\},\text{ }R_{2}:=S_{2}^{-1}%
\mathbb{Z};\text{.}
\end{eqnarray*}%
Recall from \ref{Rn} that%
\begin{eqnarray*}
Spec(R_{1}) &=&\{\{0_{R_{1}}\},(2)R_{1},(3)R_{1}\}; \\
Spec(R_{2}) &=&\{\{0_{R_{2}}\},(2)R_{2},(5)R_{2},(7)R_{2}\}.
\end{eqnarray*}%
So, $Spec(R_{1})$ is $\mathcal{V}_{2}$ and $Spec(R_{2})$ is $\mathcal{V}%
_{3}. $ By Lemma \ref{localization-AM} (1),
\begin{equation*}
Spec(R_{1}\times R_{2})\approx Spec(R_{1})\bigsqcup Spec(R_{2})
\end{equation*}%
a disjoint union $\mathcal{V}_{2}\bigsqcup \mathcal{V}_{3},$ and
\begin{equation*}
Spec(R_{1}\times R_{[2]})\approx Spec(R_{1})\bigsqcup Spec(R_{[2]})
\end{equation*}%
a disjoint union $\mathcal{V}_{2}\bigsqcup \mathcal{T}_{2}.$ It follows by
Corollary \ref{T1/4-T1/2} and Theorem \ref{t3/4t} (2) that $Spec(R_{1}\times
R_{2})$ and $Spec(R_{1}\times R_{[2]})$ are $T_{\frac{1}{2}}$ but \emph{not}
$T_{\frac{3}{4}}.\blacksquare $
\end{exs}

Notice that if the assumption "$n_{i}\geq 2$" in the first part of Theorem %
\ref{t3/4t} (1) is dropped, then the result is not necessarily true (by the
second part of the result as $\mathcal{T}_{1}=\mathcal{C}_{2}=\mathcal{V}%
_{1} $). We clarify this with a concrete example:

\begin{ex}
\label{DVR}Let $L$ be an $X$-top lattice for some $X\subseteq L\backslash
\{1\}$ with $\left\vert X\right\vert =2.$ Let $X=\{x,y\}$ for some $x\neq y$.

\textbf{Case I:}\ $x$ and $y$ are comparable, say $x<y.$ In this case, $X$
is \emph{not} $T_{\frac{3}{4}}$ as $x$ is not open regular ($int(\overline{%
\{x\}})=int(V(x))=int(X)=X\neq \{x\}$). A concrete example is $X=Spec(S),$
the prime spectrum of the local semiring $S=\{0,a,1\}\ $in Example \ref{S3}.
Moreover, this is the case for $X=Spec(D),$ where $D$ is a valuation ring
with $K.\dim (D)=1$ (e.g. a DVR).

\textbf{Case II:}$\ x$ and $y$ are \emph{not} comparable. In this case, $X$
is clearly $T_{1}$ as $X=Max(X).$ A concrete example is $X=Spec(\mathbb{F}%
_{1}\times \mathbb{F}_{2})$ for two fields $\mathbb{F}_{1}$ and $\mathbb{F}%
_{2}.\blacksquare $
\end{ex}

\begin{ex}
\label{ex-concrete}Let $(R,\mathfrak{m})$ be a \emph{local} commutative
(semi)ring with $K.\dim (R)\leq 1$ and $\left\vert Spec(R)\right\vert =n.$

\begin{enumerate}
\item If $n=1,$ then $Spec(R)=Max(X)=\{\mathfrak{m}\}\ $is discrete, whence $%
T_{2}$ (cf. Proposition \ref{dim-0}). A concrete example is $R=\mathbb{Z}%
_{p^{n}},$ where $p$ is a prime number as $Spec(\mathbb{Z}_{p^{n}})=Max(%
\mathbb{Z}_{p^{n}})=\{(p)\}.$

\item If $n=2,$ then $Spec(R)$ is $1$-dimensional and $\mathcal{T}_{1}=%
\mathcal{C}_{2}=\mathcal{V}_{1}$ whence $T_{\frac{1}{2}}$ but not $T_{\frac{3%
}{4}}$ by Theorem \ref{t3/4t} (2). Concrete examples are obtained by letting
$R$ be a valuation ring $D$ with $K.\dim (D)=1$ (e.g. a DVR; see Example \ref%
{DVR}). A $1$-dimensional local commutative ring that is \emph{not} a domain
is $R:=(\mathbb{C}[x,y]/(x^{2},xy))_{\mathfrak{m}},$ the localization of $%
\mathbb{C}[x,y)/(x^{2},xy)$ at the maximal ideal $\mathfrak{m}%
=(x,y)/(x^{2},xy).$

\item If $n\geq 3,$ then $Spec(R)$ is $1$-dimensional and $\mathcal{T}_{n-1}$
whence $T_{\frac{3}{4}}$ but not $T_{1}$ by Theorem \ref{t3/4t} (1). A
concrete example is $R:=\mathbb{F}[[x,y]]/(xy)$ (where $\mathbb{F}$ is a
field). Notice that $Spec(R)=\{\overline{(x)},$ $\overline{(y)},$ $\overline{%
(x,y)}\}$ is $\mathcal{T}_{2}$, whence $T_{\frac{3}{4}}$ but not $T_{1}$ by
Theorem \ref{t3/4t} (1).
\end{enumerate}
\end{ex}

\begin{ex}
\label{Seec-W}Consider $X=Spec(\mathbb{W}),$ the prime spectrum of the
semiring $\mathbb{W}$ and $p_{1},\cdots ,p_{n}$ ($n\geq 2$) be prime
numbers. By Corollary \ref{Y-top}, $\mathbb{W}$ is a $Y$-top semiring for
any $Y\subseteq X.$

Notice that $Y:=\{\mathbb{W}\backslash \{1\},p_{1}\mathbb{W},p_{2}\mathbb{W}%
,...,p_{n}\mathbb{W}\}$ is a $\mathcal{T}_{n}$, whence $Y$ is $T_{\frac{3}{4}%
}$ but not $T_{1}$ by Theorem \ref{t3/4t} (1).

On the other hand, $Z:=\{0,p_{1}\mathbb{W},p_{2}\mathbb{W},...,p_{n}\mathbb{W%
}\}$ is a $\mathcal{V}_{n},,$ whence $Z$ is $T_{\frac{1}{2}}$ but not $T_{%
\frac{3}{4}}$ by Theorem \ref{t3/4t} (2).$\blacksquare $
\end{ex}

An important family of semidomains was introduced by Alarcon and Anderson
\cite{AA1994} as a generalization of the rings $\mathbb{Z}_{n}=B(n,0)$ for $%
n\geq 2:$

\begin{ex}
(\cite{AA1994}, \cite[Example 1.8]{Gol1999}) Consider $B(n,i):=(B(n,i),%
\oplus ,0,\odot ,1),$ where ${B(n,i)=\{0,1,2,...,n-1\}}$ and:

\begin{enumerate}
\item $x\oplus y=x+y\ $ if $\ x+y\leq n-1$; otherwise, $x+y=u$ the unique
positive integer satisfying $i\leq u\leq n-1$ and $x+y\equiv u\ $mod\ $(n-i)$%
;

\item $x\odot y=xy\ $ if $\ xy\leq n-1$; otherwise, $xy=v$ the unique
positive integer satisfying $i\leq v\leq n-1$ and $xy\equiv v\ $mod\ $(n-i)$.%
\newline
Then $B(n,i)$ is a commutative semiring. Observe that $B(2,1)=\mathbb{B}$
and $B(n,0)=\mathbb{Z}_{n}$ for $n\geq 2.$
\end{enumerate}
\end{ex}

We recall a useful description of the prime spectra of the semidomain $%
B(n,i) $, defined in Section 1. This will be used in investigating the
separation axioms for $Spec(B(n,i)).$

\begin{thm}
\emph{(\cite[Theorem 24]{AA1994})} Let $n\geq 2,$ $1\leq i\leq n-1$ and set
\begin{equation*}
\mathfrak{m}_{n}:=\{0,2,3,...,n-1\}\text{ for }n\geq 3).
\end{equation*}

\begin{enumerate}
\item $K.\dim (B(n,i))=0$ if $i=0$ or $n=2$ and $i=1.$

\item $K.\dim (B(n,i))=1$ if $n\geq 3$ and $i=1$. In this case,%
\begin{equation*}
Spec(B(n,i)=\{0\}\cup \{pB(n,i)\text{ }|\text{ }p\text{ is a prime divisor
of }n-1\}.
\end{equation*}

\item $K.\dim (B(n,i))=1$ if $n\geq 3$ and $i=n-1$. In this case, $%
Spec(B(n,i)=\{0,\mathfrak{m}_{n}\}.$

\item $K.\dim (B(n,i))=2$ if $n\geq 4$ and $2\leq i\leq n-2$. In this case,
\begin{equation*}
Spec(B(n,i)=\{0,\mathfrak{m}_{n}\}\cup \{pB(n,i)\text{ }|\text{ }p\text{ is
a prime divisor of }n-i\}.
\end{equation*}
\end{enumerate}
\end{thm}

\begin{notatation}
For a positive integer $k\geq 2,$ we denote by $\omega (k)$ the number of
prime divisors of $k.$
\end{notatation}

\begin{cor}
\label{Bni}Let $n\geq 2$ be a positive integer, $0\leq i\leq n-1$ an integer
and set $\mathfrak{m}_{n}:=\{0,2,3,...,n-1\}$ (in case $n\geq 3$).

\begin{enumerate}
\item $Spec(\mathbb{Z}_{n}))=Max(\mathbb{Z}_{n})$ is discrete, whence $T_{2}$
\emph{(cf. Example \ref{Zn})}.

\item $Spec(\mathbb{B})=Max(\mathbb{B})=\{\{0\}\}$ is discrete, whence $%
T_{2}.$

\item If $n\geq 3,$ then $Spec(B(n,1))$ is $\mathcal{V}_{\omega (n-1)},$
whence $T_{\frac{1}{2}}$ but not $T_{\frac{3}{4}}$.

\item If $n\geq 3,$ then $Spec(B(n,n-1))$ is $\mathcal{C}_{2}=\mathcal{V}%
_{1},$ whence $T_{\frac{1}{2}}$ but not $T_{\frac{3}{4}}$.

\item If $n\geq 4$ and $2\leq i\leq n-2,$ then $Spec(B(n,i))$ is $T_{0}$ but
not $T_{\frac{1}{4}}$.
\end{enumerate}
\end{cor}

\begin{ex}
Let $n\geq 4$ and $2\leq i\leq n-2.$ Consider the semidomain $B(n,i).$ Then $%
K.\dim (B(n,i))=2$ and%
\begin{equation*}
Spec(B(n,i)=\{0,\mathfrak{m}_{n}\}\cup \{pB(n,i)\text{ }|\text{ }p\text{ is
a prime divisor of }n-i\}.
\end{equation*}%
Set $Y:=Spec(B(n,i))\backslash \{0\}.$ It follows by Corollary \ref{Y-top}
that $B(n,i)$ is a $Y$-top semiring. Notice that $Y$ is a $\mathcal{T}%
_{\omega (n-i)}$. If $n-i$ is prime, then $\omega (n-i)=1$ and $Y$ is $%
\mathcal{T}_{1}=\mathcal{C}_{2}=\mathcal{V}_{1},$ whence $T_{\frac{1}{2}}$
but not $T_{\frac{3}{4}}$ by Theorem \ref{t3/4t} (2). If $n-i$ is not prime,
then $\omega (n-i)\geq 2$ and $Y$ is $T_{\frac{3}{4}}$ but not $T_{1}$ by
Theorem \ref{t3/4t} (1).$\blacksquare $
\end{ex}

\end{document}